\documentclass[11pt,a4paper]{article}
\usepackage{mathpazo}
\usepackage{microtype}

\usepackage[a4paper,
            top=2.5cm, bottom=2.7cm,
            left=2.7cm, right=2.7cm,
            footskip=0.6cm]{geometry}

\usepackage{amsmath,amssymb,amsthm}
\usepackage{mathrsfs}
\usepackage{stmaryrd}
\usepackage{graphicx}
\usepackage{xcolor}
\usepackage{float}

\theoremstyle{plain}
\newtheorem{thm}{Theorem}[section]
\newtheorem{cor}[thm]{Corollary}
\newtheorem{lem}[thm]{Lemma}
\newtheorem{prop}[thm]{Proposition}

\newtheorem{rem}[thm]{Remark}
\newtheorem{defi}[thm]{Definition}

\def\sqr#1#2{{\vcenter{\vbox{\hrule height.#2pt
              \hbox{\vrule width.#2pt height#1pt \kern#1pt \vrule
              width.#2pt}
              \hrule height.#2pt}}}}

\def\ga{{\gamma}}

\def\ep{{\epsilon}}

\def\vep{{\varepsilon}}

\def\Sp{{\mathrm {Sp}}}

\def\R{{\Bbb R}}

\def\<{{\langle}}
\def\>{{\rangle}}
\def\no{\noindent}

\def\bs{\bigskip}

\def\dim{\hbox{\rm dim$\,$}}

\def\({\Big (}
\def\){\Big )}
\def\[{\Big[}
\def\]{\Big]}

\def\be{\begin{equation}}
\def\bel{\begin{equation}\label}
\def\ee{\end{equation}}
\def\bea{\begin{eqnarray}}
\def\eea{\end{eqnarray}}
\def\bt{\begin{theorem}}
\def\et{\end{theorem}}
\def\bc{\begin{corollary}}
\def\ec{\end{corollary}}
\def\bl{\begin{lemma}}
\def\el{\end{lemma}}
\def\bp{\begin{proposition}}
\def\ep{\end{proposition}}
\def\br{\begin{remark}}
\def\er{\end{remark}}
\def\ba{\begin{array}}
\def\ea{\end{array}}
\def\bd{\begin{definition}}
\def\ed{\end{definition}}

\makeatletter

\@addtoreset{equation}{section}
\makeatother

\begin{document}

\title{\bf  On the linear stability of elliptic relative equilibrium in planar restricted $(n+1)$-body problem }
\author{Chenshu Chu\thanks{e-mail: 202520294@mail.sdu.edu.cn}
\\
School of Mathematics, Shandong University,
Jinan, 250100, China
}
\date{}
\maketitle
\begin{abstract}
In this paper, we consider the restricted $(n+1)$-body problems, where the $n$ primaries forming an a homographic
motion generated by the planar central configuration, and the massless particle is attracted by this $n$ primaries. We are interested in the linear stability of
the elliptic relative equilibria of the massless particle under the pulsating coordinates. The main contribution of this paper first lies in the investigation of a general $(\alpha,\beta)$-system (\ref{a-b-type}). We obtain large-range spectral instability and small-range stability results valid for arbitrary eccentricity. The restricted $(n+1)$-body problem arises as a special case of this system. Accordingly, we derive several criteria for spectral instability and linear stability for this problem. As a final application, we establish some stability results of elliptic relative equilibria for massless particle under some special central configurations.

\end{abstract}

\bs

\no{\bf AMS Subject Classification:} 37J25, 70F10, 37J46, 53D12

\bs

\no{\bf Key Words:}  linear stability, central configuration, elliptic relative equilibrium, restricted $(n+1)$-body problem

\section{Introduction and main results}%Section 1.

For $N$ particles with masses $m_1,\cdots,m_N$, let $q_1,\cdots,q_N\in \mathbb{R}^2$ be the position vectors. Let
\bea U(q)=\sum_{1\leq i< j\leq N} \frac{m_im_j}{\|q_i-q_j\|},\; q=(q_{1},...,q_{N}),\nonumber\eea
be the potential function defined on the configuration space  $$
\Lambda=\{x=(x_1,\cdots,x_N)\in\mathbb{R}^{2N}\setminus\triangle:
\sum_{i=1}^Nm_ix_i=0 \},$$ where
$\triangle=\{x\in\mathbb{R}^{2N}:\exists\; i\neq j,x_i=x_j \}$ is the
collision set, $||\cdot||$ is the Euclidean distance.  Obviously, the orbits of the $N$ bodies satisfy the following Newton equation
\bea m_i\ddot{q}_i(t)=\frac{\partial
U}{\partial q_i}(q_1,...,q_N). \nonumber\label{1.2} \eea

A planar central configuration of $n$ particles with mass center in original point is formed by $n$ position vector $a=(a_1,...,a_n)\in \mathbb{R}^{2n}$ which satisfy
\bea -\mu \mathcal{M}a=\frac{\partial U(q)}{\partial q}\big{|}_{q=a}\ \ , \label{CC}\eea
for constant $\mu=U(a)/\mathcal{I}(a)>0$, where $\mathcal{M}=diag(m_{1},m_{1},m_{2},m_{2}\ldots,m_{n},m_{n})$,
$\mathcal{I}(a)=\sum m_j\|a_j\|^2$ is the moment of inertia. In this paper, we focus on the restricted $(n+1)$-body problem with $n$ primaries $m_{i}>0$, $1\leq i\leq n$ and one massless particle $m_{n+1}=0$, where
this $n$ primaries forming an a homographic motion generated by the cental configuration $a$, i.e. each particle moves on a specific Keplerian orbit while the totality of the particles move on a homographic motion. More precisely, the homographic solution generated by the cental configuration $a$ is
\bea\label{ERE}
q_{i}(t)=r(t)R(\theta(t))a_{i},\;1\leq i\leq n,
\eea
where
\bea\label{r}
r(t)=\frac{\Omega^2/\mu}{1+\mathfrak{e}\cos\theta(t)}, \ \ r^{2}\dot{\theta}(t)=\Omega,\ \ \;R(\theta)=\left(
                                                                  \begin{array}{cc}
                                                                   \cos\theta & -\sin\theta \\
                                                                    \sin\theta & \cos\theta \\
                                                                    \end{array}
                                                                    \right),
\eea
%and $\mathfrak{R}(\theta(t))=diag(R(\theta(t)),\ldots, R(\theta(t)))$, $R(\theta(t))=\left( \begin{array}{cc} \cos \theta(t) & -\sin \theta(t)\\
%\sin \theta(t) & \cos \theta(t)\end{array}\right)$,
where $\mathfrak{e}\in[0,1)$ is the eccentricity and the constant $\Omega\neq0$ is the angular momentum.
The motion of this $n$ primaries are not affected by the massless particle $q_{n+1}(t)$ and $q_{n+1}(t)$ is attracted by this $n$ primaries,
it satisfies Newton equation
\bea\label{R3B-E}
\ddot{q}_{n+1}(t)=-\sum_{i=1}^{n}m_{i}\frac{q_{n+1}-q_{i}}{||q_{n+1}-q_{i}||^3},
\eea
where $q_{i}(t)$ is given by (\ref{ERE}). For system (\ref{R3B-E}), we are interested in the linear stability of the special solutions in following form
\bea\label{qN}
q_{n+1}(t)=r(t)R(\theta(t))a_{n+1},
\eea
where $r(t)$ is given by (\ref{r}). Directed computation shows that (\ref{R3B-E}) has solutions (\ref{qN}) if and only if $a_{n+1}$ satisfies the following
algebra equation
\bea \label{RCC}
-\mu a_{n+1}=\sum_{i=1}^{n}m_{i}\frac{a_{i}-a_{n+1}}{||a_{n+1}-a_{i}||^3},\nonumber
\eea
where $a=(a_{1},...,a_{n})$ satisfies equation (\ref{CC}), this is the central configuration equation of the massless body. We
call $q_{n+1}(t)$ the \emph{Relative Equilibria} (\textbf{RE} for short) of the restricted $(n+1)$-body problem if
the eccentricity $\mathfrak{e}=0$, then this solutions are equilibria in rotating coordinates. For $\mathfrak{e}\in(0,1)$,
we call them \emph{Elliptic Relative Equilibria} (\textbf{ERE} for short), this solutions is an equilibria in pulsating coordinates.
This is a very meaningful topic, since we can image the massless particle as asteroids or satellite which affected by this $n$ primaries,
the linear stability of ERE of the restricted $(n+1)$-body problem  becomes very important.

When $n=2$, this model reduces to the classical restricted three-body problem. It is well known that the restricted three-body problem admits five relative equilibrium points, denoted by $L_{1}$ through $L_{5}$. Among them, $L_{1}$, $L_{2}$ and $L_{3}$ are collinear with the two primary bodies and are spectrally unstable, while $L_{4}$ and $L_{5}$ each form an equilateral triangle with the two primaries, respectively. These latter two are linearly stable when the mass ratio of the two primary bodies is sufficiently large.

For $n=3$, there are only two kinds of central configurations, the Lagrangian equilateral triangle central
configuration and Euler collinear central configurations. There are many works on the linear stability of the elliptic Lagrangian orbits and elliptic Euler orbits, please refer to \cite{HLS}, \cite{MSS1} , \cite{ZL} and reference therein. For the linear stability of the ERE of a fourth massless body under the influence of the homographic motion of the central configuration generated by three primary bodies, in \cite{Lea}, Leandro studied the central configurations of the restricted $4$-body problem. The
stability and instability of the circular orbits (i.e.$\mathfrak{e}=0$) of the massless body are also discussed when the configurations possess symmetry. Further, Liu and Zhou \cite{LiuZhou} obtain the linear stability of ERE of the massless body in the planar restricted $4$-body problem.

For $n\geq 4$, it is difficult to find all central configurations. There are two typical central configurations, one is the Euler-Moulton
collinear central configuration, the other is the regular $(1+\mathfrak{n})$-gon central configuration where $\mathfrak{n}$ equal masses are at the vertices of a regular $\mathfrak{n}$-gon with an additional mass at the center. The linear stability of homographic solution of Euler-Moulton central configuration and the regular $(1+\mathfrak{n})$-gon central
configuration can be found in \cite{ZL2} and \cite{HLO}, respectively. For the ERE of the massless body in
the restricted $(n+1)$-body problems, where the $n$ primaries form a regular $(1+\mathfrak{n})$-gon central configuration,
Bang and Elmabsout \cite{BE} determine the positions of relative equilibrium for the massless body and study their stability for the case $\mathfrak{e}=0$. Based on the
work in \cite{BE}, Xie, Liu and Zhou \cite{XLZ} further study the linear stability of the ERE in the restricted $(n+1)$-body problems with $\mathfrak{e}\in(0,1)$, where the $n$ primaries form an Euler-Moulton collinear central configuration or a regular $(1+\mathfrak{n})$-gon central configuration.

In this paper, first lies in the investigation of a general $(\alpha,\beta)$-system (\ref{a-b-type}). We obtain large-range spectral instability and small-range stability results valid for arbitrary eccentricity. The restricted $(n+1)$-body problem arises as a special case of this system. Accordingly, we derive several criteria for spectral instability and linear
stability for this problem. As a final application, we establish  some stability results of
elliptic relative equilibria for massless particle under some special central configurations.

\vspace{2mm}
\textbf{Notations}. For the reader's convenience, let us introduce some common notations that we shall use throughout the paper.
\begin{itemize}
\item  $\mathbb{R}$, $\mathbb{C}$, $\mathbb{N}$, $\mathbb{Z}$ and $\mathbb{U}$ denote the sets of real numbers, complex numbers, non-negative integers, and integers, as well as the unit circle, respectively;
\item $\sigma(C)$ denote the spectrum of matrix $C$; $C^{\mathcal{T}}$ denotes the transpose of the matrix $C$; $\nabla f$ denotes the gradient of function $f$ with respect to the Euclidean inner product, expressed as a column vector; $D^2 f$ denote the Hessian of $f$.
%\item
%$\nabla f$ denotes the gradient of function $f$ with respect to the Euclidean inner product, expressed as a column vector, $D^2 f$ denote the Hessian of $f$, and $C^{\mathcal{T}}$ denote the transpose matrix of $C$.
\item For two linear symmetric operators $A$ and $B$, we write $A\geq B$ if $A-B$ is positive semi-definite, and
$A>B$ if $A-B$ is positive definite.
\item  $I_j$ denotes the identity matrix on $\mathbb{R}^j$,
$\mathbb{J}_n=diag(J_2,...,J_2)_{2n\times2n}$,
$$J_{2j}=\left( \begin{array}{cccc}0_j& -I_j \\
                                  I_j& 0_j \end{array}\right),\ \ R(\theta)=\left(
                                                                                                                            \begin{array}{cc}
                                                                                                                              \cos\theta & -\sin\theta \\
                                                                                                                              \sin\theta & \cos\theta \\
                                                                                                                            \end{array}
                                                                                                                          \right).$$
Subscripts of $I$ and $\mathbb{J}$ are sometimes omitted for simplicity when clear from context.
\item $GL(\R^{2n})$ and $\text{Sym}(\R^{2n})$ denote the sets of invertible matrices and symmetric matrices in $\mathbb{R}^{2n}$, respectively.
$\Sp(2n)=\{M\in GL(\R^{2n}),  M^\mathcal{T}JM=J\} $
denotes the symplectic group.
\item Let $M_1=\left( \begin{array}{cccc}A_1& A_2 \\
A_3 & A_4 \end{array}\right)$,  $M_2=\left( \begin{array}{cccc}B_1& B_2 \\
B_3 & B_4 \end{array}\right)$ be matrices in $\Sp(2n)$ . The symplectic sum $\diamond$ is defined by
\bea  M_1\diamond M_2=\left( \begin{array}{cccccccc}A_1& 0 &A_2 & 0 \\
                                                     0 & B_1& 0  & B_2 \\
                                                    A_3& 0 &A_4 & 0 \\
                                                     0 & B_3& 0  & B_4 \end{array}\right), \nonumber   \eea
and the notation $M_1\approx M_2$ denotes that there exists a $P\in\Sp(2n)$ such that $M_1=P^{-1}M_2P$.
%We set  ${\mathcal R}(\theta_{1})=\left( \begin{array}{cc} \cos \theta_{1} & -\sin \theta_{1}\\
%\sin \theta_{1} & \cos \theta_{1}\end{array}\right)\in\Sp(2)$.
%Then $M\in\Sp(2n)$ is linearly stable if and only if $$M\approx {R}(\theta_1)\diamond\cdots\diamond {R}(\theta_n), \quad \theta_i\in [0,2\pi).$$
\end{itemize}
In what follows, we shall introduce several definitions and present our main theorem.
\begin{defi}
A symplectic matrix $W\in\Sp(2k)$ is called spectrally stable if all eigenvalues of $W$ belong to the
unit circle $\mathbb{U}$ of the complex plane. $W$ is called linearly stable if it is spectrally stable and semi-simple.
$W$ is called elliptic-hyperbolic if some of its eigenvalues lie on the unit circle $\mathbb{U}$ and some lie outside the unit circle.
While $W$ is called hyperbolic if no eigenvalues of $W$ are on $\mathbb{U}$.
\end{defi}
Let $T>0$, for a $T$-periodic solution $x(t)$ of the following Hamiltonian system
\bea
\dot{x}(t)=J\nabla H(x(t)),\; x(0)=x(T).\label{Ham1}
\eea
Let $\ga$ be the fundamental solution of the linear Hamiltonian system corresponding to (\ref{Ham1}),
that is $$\dot{\ga}=JD^{2}H(x(t))\ga,\; \ga(0)=I.$$ then $\ga(T)$ is a symplectic matrix.
\begin{defi}
The periodic solution $x(t)$ is called hyperbolic
(spectrally stable, linearly stable, elliptic-hyperbolic, resp.) if the monodromy matrix $\ga(T)$
is hyperbolic (spectrally stable, linearly stable, elliptic-hyperbolic, resp.).
\end{defi}
In this paper , we first consider the following $(\alpha,\beta)$-type Hamiltonian system
\bea
\xi'_{\alpha,\beta,\mathfrak{e}}(\theta)=J_{4}\mathfrak{B}_{\alpha,\beta,\mathfrak{e}}(\theta)\xi_{\alpha,\beta,\mathfrak{e}}(\theta),\; \xi_{\alpha,\beta,\mathfrak{e}}(0)=I_{4}.\label{a-b-type}
\eea
with
\bea \mathfrak{B}_{\alpha,\beta,\mathfrak{e}}(\theta)=\left( \begin{array}{cccc} I_{2} & -J_{2} \\
J_{2} & I_{2}-\frac{\mathfrak{R}_{\alpha, \beta}}{1+\mathfrak{e}\cos\theta}
\end{array}\right),\quad \theta\in[0,2\pi],
\label{a-b-type-B}
\eea
where
\bea
\mathfrak{R}_{\alpha, \beta}=(1+\alpha)I_{2}+\beta\mathcal{N},\ \ \mathcal{N}=\left(
                                                             \begin{array}{cc}
                                                               1 & 0 \\
                                                               0 & -1 \\
                                                             \end{array}
                                                           \right),\nonumber
\eea
and $\alpha\geq0$, $\beta\geq0$ are the parameters, $\mathfrak{e}$ is the eccentricity. Then we have following spectral instability and linear stability results,
\begin{thm}\label{stable-thm}
%\vskip 8pt
(I)\; For $\beta>\alpha+1$, $\alpha>0$, $\xi_{\alpha,\beta,\mathfrak{e}}(2\pi)\approx D(\lambda)\diamond M$ for any $\mathfrak{e}\in[0,1)$, where $D(\lambda) = \begin{pmatrix} \lambda & 0 \\ 0 & \lambda^{-1} \end{pmatrix}$, $1\neq\lambda>0$ and $M\in\Sp(2)$;
\vskip 5pt
(II)\; For any fixed $\mathfrak{e}\in[0,1)$, let $k(\mathfrak{e})=1+(1+3\mathfrak{e}/2)^{-1}$, when $\alpha>0, 0\leq\beta<\min\{k(\mathfrak{e})\alpha,\alpha+1\}$, $\xi_{\alpha,\beta,\mathfrak{e}}(2\pi)$ is hyperbolic. In particular,
for $0\leq\beta<\min\{\frac{7}{5}\alpha,\alpha+1\}$, $\xi_{\alpha,\beta,\mathfrak{e}}(2\pi)$ is hyperbolic for any $\mathfrak{e}\in[0,1)$.
\end{thm}
\begin{thm}\label{stable-thm2}
Fixing $\mathfrak{e}\in[0,1)$, for $\alpha>0, 0\leq\beta<\min\{3\alpha,\alpha+1\}$, when $(\alpha,\beta)$ is sufficiently close to $(1/2,3/2)$, $\xi_{\alpha,\beta,\mathfrak{e}}(2\pi)$ is linear stable and $\xi_{\alpha,\beta,\mathfrak{e}}(2\pi)\approx R(2\pi-\vartheta_{1})\diamond R(2\pi-\vartheta_{2})$ with some $\vartheta_{1},\vartheta_{2}\in(0,\pi)$.
\end{thm}
\begin{figure}[H]
  \centering
   \includegraphics[height=0.52\textwidth,width=0.66\textwidth]{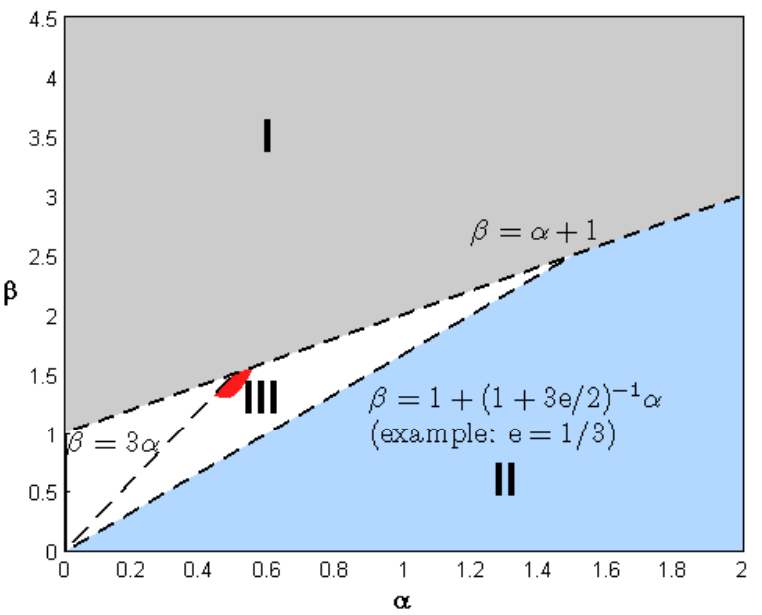}
     \caption{The spectral instability region $\textbf{I}$, the hyperbolic region $\textbf{II} (\text{example}:\mathfrak{e}=1/3)$ and the linear stable region $\textbf{III}$.}
\end{figure}
In order to study the linear stability of the ERE in the planar restricted $(n+1)$-body problem, we first need the following theorem.
\begin{thm}(See Theorem 1.1 of Xie, Liu and Zhou in \cite{XLZ})\label{red-ERE}
For the elliptic relative equilibria  $q_{n+1}(t)=r(t)R(\theta(t))a_{n+1}$ of the planar restricted $(n+1)$-body problem which produced by the planar central configuration $a=(a_{1},...,a_{n})$ of $n$ primaries, after some linear symplectic transformations, the fundamental solution of the ERE satisfies the following equation,
\bea
\zeta'_{\mathfrak{e}}(\theta)=J_{4}\mathfrak{B}_{\hat{\alpha},\hat{\beta},\mathfrak{e}}(\theta)\zeta_{\mathfrak{e}}(\theta),\; \zeta_{\mathfrak{e}}(0)=I_{4}.\label{EREH}
\eea
%with
%\bea B(\theta)=\left( \begin{array}{cccc} I_{2} & -J_{2} \\
%J_{2} & I_{2}-\frac{(1+\alpha) I_{2}+\beta\mathcal{N}}{1+\mathfrak{e}\cos\theta}
%\end{array}\right),\quad \theta\in[0,2\pi],\ \ \mathcal{N}=\left(
%                                                             \begin{array}{cc}
%                                                               1 & 0 \\
%                                                               0 & -1 \\
%                                                             \end{array}
%                                                           \right),
% \label{msf} \eea
with $\hat{\alpha}=\frac{\lambda_{\min}+\lambda_{\max}}{2}>0$, $\hat{\beta}=\frac{\lambda_{\max}-\lambda_{\min}}{2}\geq0$, $\lambda_{\min}$ and $\lambda_{\max}$ are the minimal and maximal eigenvalues of
$$
\mu^{-1}D^2\mathcal{U}(a_{n+1})=-\frac{1}{\mu}\sum_{i=1}^{n}\frac{m_{i}}{||a_{i}-a_{n+1}||^3}I_{2}+\frac{3}{\mu}\sum_{i=1}^{n}\frac{m_{i}(a_{i}-a_{n+1})(a_{i}-a_{n+1})^{\mathcal{T}}}{||a_{i}-a_{n+1}||^5}
$$
where $\mu=\frac{U(a)}{\mathcal{I}(a)}, \mathcal{U}(z)=\sum_{i=1}^{n}\frac{m_{i}}{||a_{i}-z||}$.
\end{thm}
\begin{rem}
In paper \cite{XLZ}, Xie, Liu and Zhou first apply the symplectic reduction method of Meyer and Schmidt \cite{MS} to obtain
the linear Hamiltonian system for the massless body, i.e
\bea\label{ZForm}
\zeta'_{\mathfrak{e}}(\theta)=J_{4}B(\theta)\zeta_{\mathfrak{e}}(\theta),\;
B(\theta)=\left( \begin{array}{cccc} I_{2} & -J_{2} \\
J_{2} & I_{2}-\frac{I_{2}+\mu^{-1}D^2\mathcal{U}(a_{n+1})}{1+\mathfrak{e}\cos\theta}
\end{array}\right),
\eea
and they  analysis the linear stability of two special cases:  the $n$ primaries form an Euler-Moulton collinear central configuration or a regular $(1+\mathfrak{n})$-gon central configuration.

In this paper, for the reader's convenience, we  re-derive the linear Hamiltonian system (\ref{EREH}) for the massless body, by using the linear symplectic transformations directly. Furthermore, we introduce a more general $(\alpha, \beta)$-type system (\ref{a-b-type-B}), of which (\ref{ZForm}) is merely a special case. We carry out some spectral instability and linear stability analysis on (\ref{a-b-type-B}) and finally apply them to the specific instance (\ref{ZForm}).
A special case of the $(\alpha,\beta)$-type system is first studied by Hu, Long and Ou \cite{HLO} in their
study of the linear stability of the regular $(1+\mathfrak{n})$-gon system. %In this paper, we obtain a more complete analysis for
%this $(\alpha, \beta)$-type system.
\end{rem}
As an application of Theorem \ref{stable-thm} and Theorem \ref{stable-thm2}, we have following spectral instability and linear stability results in the planar restricted $(n+1)$-body problem
\begin{thm}\label{unstable-ERE}
Let $\lambda_{\max}$ and $\lambda_{\min}$ be the largest and smallest eigenvalues of $\mu^{-1}D^2\mathcal{U}(a_{n+1})$, respectively. Then, for any eccentricity $\mathfrak{e}\in[0,1)$, the elliptic relative equilibria $q_{n+1}(t)=r(t)R(\theta(t))a_{n+1}$ of the planar restricted $(n+1)$-body problem is spectral instability for $\lambda_{\min}<-1$, and hyperbolic for $\lambda_{\min}>\max\{-\frac{\lambda_{\max}}{3+3\mathfrak{e}},-1\}$.
\end{thm}
\begin{thm}\label{stable-ERE}
Let $\lambda_{\max}$ and $\lambda_{\min}$ be the largest and smallest eigenvalues of $\mu^{-1}D^2\mathcal{U}(a_{n+1})$,respectively. Fixing eccentricity $\mathfrak{e}\in[0,1)$, when
$\lambda_{\min}>\max\{-\frac{\lambda_{\max}}{2},-1\}$ with $(\lambda_{\max},\lambda_{\min})$ sufficiently close to $(2,-1)$, then
the elliptic relative equilibria $q_{n+1}(t)=r(t)R(\theta(t))a_{n+1}$ of the planar restricted $(n+1)$-body problem is linear stability.
\end{thm}
In particular, as an application of Theorem \ref{unstable-ERE}, we consider $\mathfrak{n}$ particles, each with mass $m$, positioned at the vertices of a regular polygon, rotating around $m_{0}$ at its center. Here note that
$n=\mathfrak{n}+1, M=\mathfrak{n}m$, In \cite{BE}, Bang and Elmabsout determine the positions of relative equilibrium for the massless body lying (up to rotations of angle $2\pi/\mathfrak{n}$)
on the semi-axes $\vartheta=0$ and $\vartheta=\pi/\mathfrak{n}$, (see Fig. 2):
\begin{description}
  \item[(i)] On the semi-axis $\vartheta=0$, there are two solutions $\rho=S_{1}>1$, and $\rho=S_{2}\in(0,1);$
  \item[(ii)] On the semi-axis $\vartheta=\pi/\mathfrak{n}$, there is one solution $\rho=S_{3}>1$.
\end{description}
Here, $S_{1}, S_{2}$ and $S_{3}$ also denote the positions of the massless body, respectively. In the same paper, Bang and Elmabsout have proved that for $\mathfrak{e}=0$,
$S_{1}$ and $S_{2}$ are spectrally unstable for every $m_{0}>0$, while $S_{3}$ is linearly stable for $n\geq7$ if $m_{0}/M$
is sufficiently large.

Recently, based on the symplectic reduction and the Maslov-type index theory, Xie, Liu and Zhou \cite{XLZ} study the elliptic relative equilibrium with eccentricity $\mathfrak{e}\in(0,1)$ and further prove that
\begin{thm}\cite[Theorem 1.4]{XLZ}\label{XLZ-TH1}
If the massless body lies at $S_1$ or $S_{2}$, there exists an $\mathfrak{e}_{0}\in(0,1)$ such that the elliptic relative equilibrium is linearly unstable
when $m_{0}/M$ is sufficiently large and the eccentricity $\mathfrak{e}\in[0,\mathfrak{e}_{0}]$.
\end{thm}
\begin{thm}\cite[Theorem 1.5]{XLZ}\label{XLZ-TH2}
if the massless body lies at $S_{3}$, the elliptic relative equilibrium is linear stability when $m_{0}/M$
is sufficiently large.
\end{thm}
\begin{figure}[H]
  \centering
   \includegraphics[height=0.4\textwidth,width=0.6\textwidth]{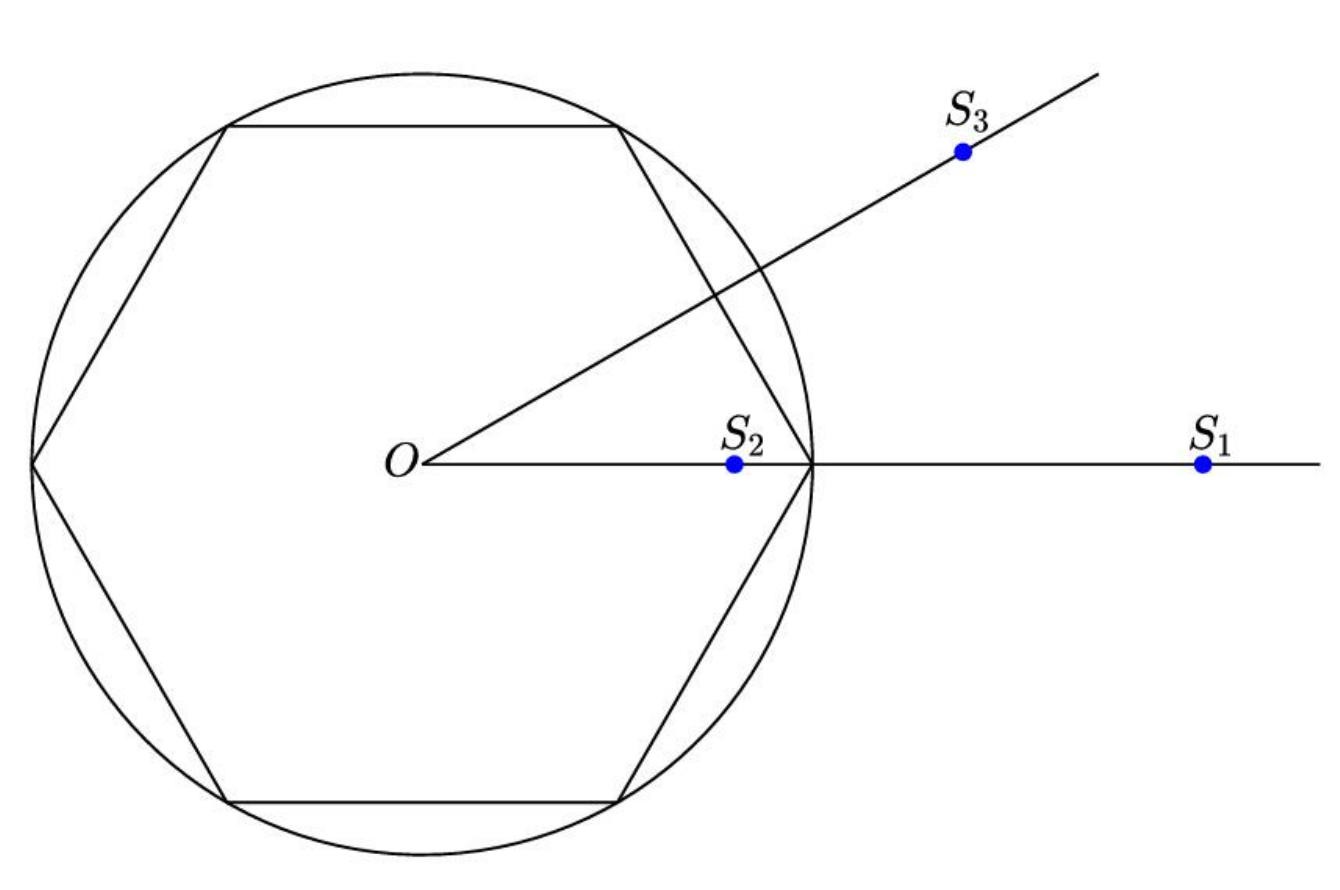}
     \caption{Positions of relative equilibrium for the massless body in the case $m_{0}>0$.}
\end{figure}
In this paper, based on our Theorem \ref{stable-ERE}, we re-prove Theorem \ref{XLZ-TH2}. Moreover, base on our Theorem \ref{unstable-ERE}, we improved the result of Theorem \ref{XLZ-TH1} and further obtained
\begin{thm}\label{1+n-gon}
If the massless body lies at $S_1$ or $S_{2}$, the elliptic relative equilibrium of the massless body is spectral unstable for any $\mathfrak{e}\in[0,1)$.
\end{thm}

This paper is organized as follows. In Section 2, we first study the linear stability of the $(\alpha,\beta)$-system, based on the Maslov-type index theory, we obtain our main Theorem \ref{stable-thm} and \ref{stable-thm2}, which provides a large spectral-instability region and a small stability region, both valid for arbitrary eccentricity
$\mathfrak{e}\in[0,1)$.
In Section 3, we study the linear stability of the planar restricted $(n+1)$-body problem, first we utilize a linear symplectic transformation to re-derive the linear system for the ERE of the restricted $(n+1)$-body problem given in Xie-Liu-Zhou \cite{XLZ}. Subsequently, we convert this linear system into the special form of the $(\alpha,\beta)$-system, and then employ Theorem \ref{stable-thm} and \ref{stable-thm2} to derive the stability criterion (i.e.Theorem \ref{unstable-ERE} and Theorem \ref{stable-ERE}) for the ERE in the restricted $(n+1)$-body problem.
In Section 4, based on Theorem \ref{unstable-ERE} and Theorem \ref{stable-ERE}, we study the linear stability of the ERE in the restricted $(n+1)$-body problems with any $\mathfrak{e}\in[0,1)$, where the $n$ primaries form in some specific central configurations. In the appendix to Section 5, we briefly review the Maslov-type index and its properties, as well as its relationship with the relative Morse index.
\section{The spectral instability of the $(\alpha,\beta)$-system}
For the $(\alpha,\beta)$-system with $\alpha\geq0, \beta\geq0$, we first have the following simple proposition.
\begin{prop}\label{spec-e=0}
For the eccentricity $\mathfrak{e}=0$, the monodromy matrix $\xi_{\alpha,\beta,0}(2\pi)$ has following properties:

(i)\; The eigenvalues of $\xi_{\alpha,\beta,0}(2\pi)$ are
\bea
\sigma(\xi_{\alpha,\beta,0}(2\pi))=\left\{e^{\pm2\pi\sqrt{\alpha-1+\sqrt{\beta^2-4\alpha}}},\;
e^{\pm2\pi\sqrt{\alpha-1-\sqrt{\beta^2-4\alpha}}}\right\}.\nonumber
\eea

(ii)\; $\xi_{\alpha,\beta,0}(2\pi)$ is elliptic if and only if
$$
2\sqrt{\alpha}\leq\beta\leq\alpha+1,\; 0\leq\alpha\leq1.
$$

(iii)\; $\xi_{\alpha,\beta,0}(2\pi)$ is elliptic-hyperbolic if and only if
$$
\alpha+1<\beta\;\;\text{or}\;\;(\alpha+1=\beta, \alpha>1).
$$

(iv)\; $\xi_{\alpha,\beta,0}(2\pi)$ is hyperbolic if and only if
$$
\beta<2\sqrt{\alpha}\;\; \text{or}\;\; (2\sqrt{\alpha}\leq\beta<\alpha+1, \alpha>1).
$$
\begin{figure}[H]
  \centering
   \includegraphics[height=0.5\textwidth,width=0.68\textwidth]{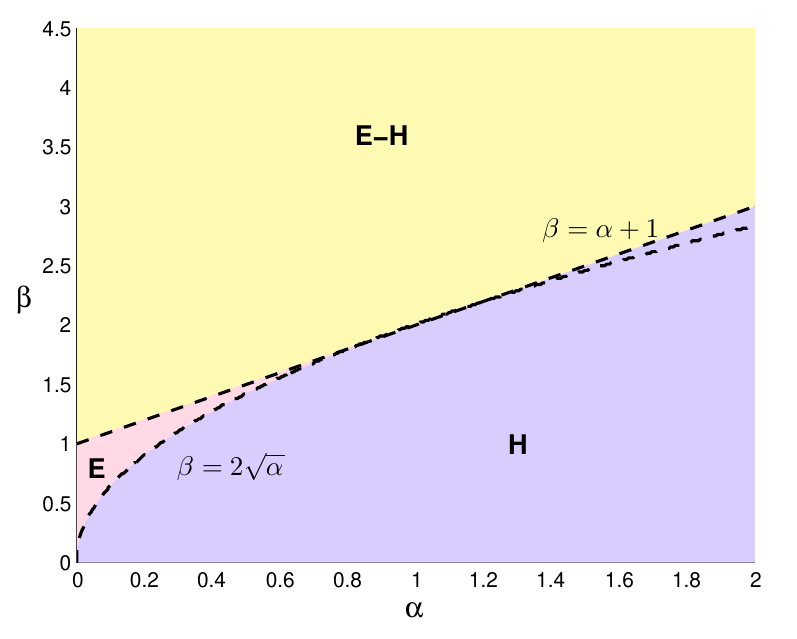}
     \caption{The elliptic (\textbf{E}), elliptic-hyperbolic (\textbf{E-H}) and hyperbolic (\textbf{H}) regions for $\mathfrak{e}=0$.}
\end{figure}
\end{prop}
\begin{proof}
The reader can easily obtain this by direct calculation, so we omit the proof details.
\end{proof}
Next, we consider the linear operator,
\bea\label{ope-est}
\mathcal{A}(R,\mathfrak{e}):=-\frac{d^2}{d\theta^2}I_{2k}-2\mathbb{J}_{k}\frac{d}{d\theta}+\frac{\mathfrak{R}}{1+\mathfrak{e}\cos\theta},\nonumber
\eea
where $\mathfrak{R}$ is a $2k\times2k$ symmetry matrix. Let
$$ \bar{D}_{2k}(\omega,2\pi)=\{y\in W^{2,2}([0,2\pi],\mathbb{C}^{2k})\;|\;y(2\pi)=\omega y(0),\dot{y}(2\pi)=\omega \dot{y}(0)\}, $$
where $\omega\in\mathbb{U}$ and $W^{2,2}([0,2\pi])$ is the usual Sobolev space.
Then $\mathcal{A}(\mathfrak{R},\mathfrak{e})$ is a self-adjoint operator in $L^2([0,2\pi], \mathbb{C}^{2k})$ with domain
$\bar{D}_{2k}(\omega,2\pi)$. We define the $\omega$-Morse index $\phi_{\omega}(\mathcal{A}(\mathfrak{R},e))$ to be the total number of negative
eigenvalues of $\mathcal{A}(\mathfrak{R},\mathfrak{e})$, and define
$\nu_{\omega}(\mathcal{A}(\mathfrak{R},\mathfrak{e}))=\dim\ker(\mathcal{A}(\mathfrak{R},\mathfrak{e}))$. In particular, when $\mathfrak{R}=\mathfrak{R}_{\alpha,\beta}$, the operator
$\mathcal{A}(\mathfrak{R},\mathfrak{e})$ corresponding the linear system (\ref{a-b-type}). Then, we have following useful propositions and lemmas,
\begin{prop}\label{rela-Morse-Maslov} (See Long \cite{Lon4} p.172).
The $\omega$-Morse index $\phi_{\omega}(\mathcal{A}(\mathfrak{R}_{\alpha,\beta},e))$ and nullity $\nu_{\omega}(\mathcal{A}(\mathfrak{R}_{\alpha,\beta},e))$ are equal
to the $\omega$-Maslov-type index $i_{\omega}(\xi_{\alpha,\beta,\mathfrak{e}})$ and nullity $\nu_{\omega}(\xi_{\alpha,\beta,\mathfrak{e}})$ respectively, that is,
for any $\omega\in \mathbb{U}$, we have
\bea \phi_{\omega}(\mathcal{A}(\mathfrak{R}_{\alpha,\beta},\mathfrak{e}))=i_{\omega}(\xi_{\alpha,\beta,\mathfrak{e}}),\ \ \nu_{\omega}(\mathcal{A}(\mathfrak{R}_{\alpha,\beta},e))=\nu_{\omega}(\xi_{\alpha,\beta,\mathfrak{e}}).\nonumber \label{eq}\eea
where $\xi_{\alpha,\beta,\mathfrak{e}}$ is given by (\ref{a-b-type}).
\end{prop}
\begin{lem}\cite[Lemma 5.5]{HO}\label{lem-pos}
The operator $\mathcal{A}(I_{2k},\mathfrak{e})$ is positive definite on its domain $\bar{D}_{2k}(\omega,2\pi)$ for any $\omega\in\mathbb{U}\setminus\{1\}$.
\end{lem}
\begin{lem}\cite[Lemma 4.1]{ZL}\label{lem-pos2}
The operator $\mathcal{A}(I_{2k},\mathfrak{e})$ is positive semi-definite on its domain $\bar{D}_{2k}(1,2\pi)$ with $\ker\mathcal{A}(I_{2k},\mathfrak{e})=\oplus_{l=1}^{k}\ker\mathcal{A}(I_{2},\mathfrak{e})$, where
$$
\ker\mathcal{A}(I_{2},\mathfrak{e})=\text{\emph{Span}}_{\mathbb{C}}\left\{\left(
                                                     \begin{array}{c}
                                                       (1+\mathfrak{e}\cos\theta)\cos\theta \\
                                                       -(1+\mathfrak{e}\cos\theta)\sin\theta \\
                                                     \end{array}
                                                   \right),\;\left(
                                                     \begin{array}{c}
                                                       (1+\mathfrak{e}\cos\theta)\sin\theta \\
                                                       (1+\mathfrak{e}\cos\theta)\cos\theta \\
                                                     \end{array}
                                                   \right)\right\}.
$$
\end{lem}
%\textbf{Proof of Theorem \ref{stable-thm}-(II)}):
%When $\alpha>0$ and $0\leq\beta<\alpha$, $\mathcal{A}(\mathfrak{R}_{\alpha,\beta},\mathfrak{e})$ is strictly positive in $\bar{D}_{2}(\omega,2\pi)$ for any $\omega\in\mathbb{U}$ and
%$\xi_{\alpha,\beta,\mathfrak{e}}(2\pi)$ is hyperbolic.

%\begin{proof}[\textbf{Proof of Theorem \ref{stable-thm}-(II)}]
%When $\alpha>0$ and $0\leq\beta<\alpha$, we have $\mathfrak{R}_{\alpha,\beta}\geq(1+\alpha-\beta)I_{2}$, then Lemma \ref{lem-pos} implies that $\mathcal{A}(\mathfrak{R}_{\alpha,\beta},\mathfrak{e})$ is strictly positive in $\bar{D}_{2}(\omega,2\pi)$ for any $\omega\in\mathbb{U}$. This further tells us that $\nu_{\omega}(\xi_{\alpha,\beta,\mathfrak{e}}(2\pi))=0$ for all $\omega\in\mathbb{U}$, hence $\sigma(\xi_{\alpha,\beta,\mathfrak{e}}(2\pi))\cap \mathbb{U}={\o}$, then
%$\xi_{\alpha,\beta,\mathfrak{e}}(2\pi)$ is hyperbolic. This completes the proof.
%\end{proof}
\begin{prop}\label{monot}
For any fixed $\mathfrak{e}\in[0,1)$, $\alpha>0$, $\beta>0$, we have
$$
\begin{aligned}
&i_{\omega}(\xi_{\alpha,\beta+\epsilon,\mathfrak{e}})-i_{\omega}(\xi_{\alpha,\beta,\mathfrak{e}})\geq\nu_{\omega}(\xi_{\alpha,\beta,\mathfrak{e}}),\ \ \nu_{\omega}(\xi_{\alpha,\beta,\mathfrak{e}})=\lim_{\epsilon\to0^{+}}i_{\omega}(\xi_{\alpha,\beta+\epsilon,\mathfrak{e}})-i_{\omega}(\xi_{\alpha,\beta,\mathfrak{e}}),\\
&i_{\omega}(\xi_{\alpha-\epsilon,\beta,\mathfrak{e}})-i_{\omega}(\xi_{\alpha,\beta,\mathfrak{e}})\geq\nu_{\omega}(\xi_{\alpha,\beta,\mathfrak{e}}),\ \ \nu_{\omega}(\xi_{\alpha,\beta,\mathfrak{e}})=\lim_{\epsilon\to0^{+}}i_{\omega}(\xi_{\alpha-\epsilon,\beta,\mathfrak{e}})-i_{\omega}(\xi_{\alpha,\beta,\mathfrak{e}}).
\end{aligned}
$$
\end{prop}
\begin{proof}
First, we have
$$
\phi_{\omega}(\mathcal{A}(\mathfrak{R}_{\alpha,\beta},\mathfrak{e}))=\phi_{\omega}(\beta^{-1}\mathcal{A}(\mathfrak{R}_{\alpha,\beta},\mathfrak{e})),\;
\nu_{\omega}(\mathcal{A}(\mathfrak{R}_{\alpha,\beta},\mathfrak{e}))=\nu_{\omega}(\beta^{-1}\mathcal{A}(\mathfrak{R}_{\alpha,\beta},\mathfrak{e})),
$$
direct computations show that
$$
\beta^{-1}\mathcal{A}(\mathfrak{R}_{\alpha,\beta},\mathfrak{e})=\beta^{-1}\mathcal{A}(\mathfrak{R}_{\alpha,0},\mathfrak{e})+\frac{\mathcal{N}}{1+\mathfrak{e}\cos\theta},
$$
hence for $\epsilon>0$
$$
(\beta+\epsilon)^{-1}\mathcal{A}(\mathfrak{R}_{\alpha,\beta+\epsilon},\mathfrak{e})-\beta^{-1}\mathcal{A}(\mathfrak{R}_{\alpha,\beta},\mathfrak{e})
=-\epsilon\beta^{-1}(\beta+\epsilon)^{-1}\mathcal{A}(\mathfrak{R}_{\alpha,0},\mathfrak{e}).
$$
From Lemma \ref{lem-pos} and Lemma \ref{lem-pos2}, we know that $\mathcal{A}(\mathfrak{R}_{\alpha,0},\mathfrak{e})>0$ is  positive definite in $\bar{D}_{2}(\omega,2\pi)$ for any $\omega\in\mathbb{U}$, hence
$$
\phi_{\omega}((\beta+\epsilon)^{-1}\mathcal{A}(\mathfrak{R}_{\alpha,\beta+\epsilon},\mathfrak{e}))-\phi_{\omega}(\beta^{-1}\mathcal{A}(\mathfrak{R}_{\alpha,\beta},\mathfrak{e}))
\geq\nu_{\omega}(\beta^{-1}\mathcal{A}(\mathfrak{R}_{\alpha,\beta},\mathfrak{e})).
$$
This implies
\bea
\phi_{\omega}(\mathcal{A}(\mathfrak{R}_{\alpha,\beta+\epsilon},\mathfrak{e}))-\phi_{\omega}(\mathcal{A}(\mathfrak{R}_{\alpha,\beta},\mathfrak{e}))
\geq\nu_{\omega}(\mathcal{A}(\mathfrak{R}_{\alpha,\beta},\mathfrak{e})),\label{index-equ1}
\eea
or equivalent to
$$
i_{\omega}(\xi_{\alpha,\beta+\epsilon,\mathfrak{e}})-i_{\omega}(\xi_{\alpha,\beta,\mathfrak{e}})\geq\nu_{\omega}(\xi_{\alpha,\beta,\mathfrak{e}}).
$$
On the other hand, by the Proposition \ref{per-index} of the Maslov-type index, we always have
$$
\lim_{\epsilon\to0^{+}}i_{\omega}(\xi_{\alpha,\beta+\epsilon,\mathfrak{e}})-i_{\omega}(\xi_{\alpha,\beta,\mathfrak{e}})\leq\nu_{\omega}(\xi_{\alpha,\beta,\mathfrak{e}}).
$$
This completes the proof of the first equality.

For the second equality, one can see that
$$
\mathcal{A}(\mathfrak{R}_{\alpha-\epsilon,\beta},\mathfrak{e})=\mathcal{A}(\mathfrak{R}_{\alpha,\beta},\mathfrak{e})-\frac{\epsilon I_{2}}{1+\mathfrak{e}\cos\theta},
$$
hence
$$
\phi_{\omega}(\mathcal{A}(\mathfrak{R}_{\alpha-\epsilon,\beta},\mathfrak{e}))-\phi_{\omega}(\mathcal{A}(\mathfrak{R}_{\alpha,\beta},\mathfrak{e}))
\geq\nu_{\omega}(\mathcal{A}(\mathfrak{R}_{\alpha,\beta},\mathfrak{e})),
$$
or equivalent to
$$
i_{\omega}(\xi_{\alpha-\epsilon,\beta,\mathfrak{e}})-i_{\omega}(\xi_{\alpha,\beta,\mathfrak{e}})\geq\nu_{\omega}(\xi_{\alpha,\beta,\mathfrak{e}}).
$$
On the other hand, by the Proposition \ref{per-index} of the Maslov-type index, we always have
$$
\lim_{\epsilon\to0^{+}}i_{\omega}(\xi_{\alpha-\epsilon,\beta,\mathfrak{e}})-i_{\omega}(\xi_{\alpha,\beta,\mathfrak{e}})\leq\nu_{\omega}(\xi_{\alpha,\beta,\mathfrak{e}}).
$$
This completes the proof of the second equality.
\end{proof}
\begin{prop}\label{monot2}
For any fixed $\mathfrak{e}\in[0,1)$, $k\in \mathbb{R}$, $\alpha>0$, we have
$$
\begin{aligned}
i_{\omega}(\xi_{\alpha+\epsilon,k(\alpha+\epsilon),\mathfrak{e}})-
i_{\omega}(\xi_{\alpha,k\alpha,\mathfrak{e}})\geq \nu_{\omega}(\xi_{\alpha,k\alpha,\mathfrak{e}}),\ \
\nu_{\omega}(\xi_{\alpha,k\alpha,\mathfrak{e}})=\lim_{\epsilon\to0^{+}}i_{\omega}(\xi_{\alpha+\epsilon,k(\alpha+\epsilon),\mathfrak{e}})-i_{\omega}(\xi_{\alpha,k\alpha,\mathfrak{e}}).
\end{aligned}
$$
\end{prop}
\begin{proof}
Direct computation shows that
\bea\label{equal3}
\alpha^{-1}\mathcal{A}(\mathfrak{R}_{\alpha,k\alpha},\mathfrak{e})=\alpha^{-1}\mathcal{A}(I_{2},\mathfrak{e})+\frac{1}{1+\mathfrak{e}\cos\theta}\left(\begin{array}{cc}
1+k & 0 \\
0 & 1-k \\
\end{array}
\right),
\eea
then for any $\epsilon>0$,
\bea\label{equal4}
(\alpha+\epsilon)^{-1}\mathcal{A}(\mathfrak{R}_{\alpha+\epsilon,k(\alpha+\epsilon)},\mathfrak{e})-\alpha^{-1}\mathcal{A}(\mathfrak{R}_{\alpha,k\alpha},\mathfrak{e})=-
\epsilon\alpha^{-1}(\alpha+\epsilon)^{-1}\mathcal{A}(I_{2},\mathfrak{e}).
\eea
First, from Lemma \ref{lem-pos}, we know the operator $\mathcal{A}(I_{2},\mathfrak{e})$ is a positive operator on its domain $\bar{D}_{2}(\omega,2\pi)$ for any $\omega\in\mathbb{U}\setminus\{1\}$, hence $(\ref{equal4})$ implies
\bea\label{index equ3}
\phi_{\omega}(\mathcal{A}(\mathfrak{R}_{\alpha+\epsilon,k(\alpha+\epsilon)},\mathfrak{e}),\mathfrak{e}))-\phi_{\omega}(\mathcal{A}(\mathfrak{R}_{\alpha,k\alpha},\mathfrak{e}))
\geq\nu_{\omega}(\mathcal{A}(\mathfrak{R}_{\alpha,k\alpha},\mathfrak{e})),\;\omega\in\mathbb{U}\setminus\{1\}.
\eea
Second , Lemma \ref{lem-pos2} shows that $\mathcal{A}(I_{2},\mathfrak{e})$ is a positive operator on its domain $\bar{D}_{2}(1,2\pi)$ and
$$
\ker\mathcal{A}(I_{2},\mathfrak{e})=\text{\emph{Span}}_{\mathbb{C}}\left\{\left(
                                                     \begin{array}{c}
                                                       (1+\mathfrak{e}\cos\theta)\cos\theta \\
                                                       -(1+\mathfrak{e}\cos\theta)\sin\theta \\
                                                     \end{array}
                                                   \right),\;\left(
                                                     \begin{array}{c}
                                                       (1+\mathfrak{e}\cos\theta)\sin\theta \\
                                                       (1+\mathfrak{e}\cos\theta)\cos\theta \\
                                                     \end{array}
                                                   \right)\right\}.
$$
then from equality (\ref{equal3}), if $x\in\text{Ker}(\mathcal{A}(\mathfrak{R}_{\alpha,k\alpha},\mathfrak{e}))$, we must have $x\notin\text{Ker}(\mathcal{A}(I_{2},\mathfrak{e}))$, hence
$(\mathcal{A}(I_{2},\mathfrak{e})x,x)_{L^{2}}>0$. This combine with inequality (\ref{equal4}), we obtain
\bea\label{index equ4}
\phi_{1}(\mathcal{A}(\mathfrak{R}_{\alpha+\epsilon,k(\alpha+\epsilon)},\mathfrak{e}),\mathfrak{e}))-\phi_{1}(\mathcal{A}(\mathfrak{R}_{\alpha,k\alpha},\mathfrak{e}))
\geq\nu_{1}(\mathcal{A}(\mathfrak{R}_{\alpha,k\alpha},\mathfrak{e})).
\eea
Therefore, (\ref{index equ3}) and (\ref{index equ4}) combine with Proposition \ref{rela-Morse-Maslov} implies for any $\omega\in\mathbb{U}$, we have
$$
i_{\omega}(\xi_{\alpha+\epsilon,k(\alpha+\epsilon),\mathfrak{e}})-i_{\omega}(\xi_{\alpha,k\alpha,\mathfrak{e}})\geq\nu_{\omega}(\xi_{\alpha,k\alpha,\mathfrak{e}}).
$$
On the other hand, by the Proposition \ref{per-index} of the Maslov-type index, we always have
$$
\lim_{\epsilon\to0^{+}}i_{\omega}(\xi_{\alpha+\epsilon,k(\alpha+\epsilon),\mathfrak{e}})-i_{\omega}(\xi_{\alpha,k\alpha,\mathfrak{e}})\leq\nu_{\omega}(\xi_{\alpha,k\alpha,\mathfrak{e}}).
$$
This completes the proof.
\end{proof}
\begin{cor}\label{isolate}
(1) For any fixed $\mathfrak{e}\in[0,1)$ and $\alpha>0$, the point $\beta>0$ satisfying $\nu_{\omega}(\xi_{\alpha,\beta,\mathfrak{e}})\neq0$ is isolated;
(2) For any fixed $\mathfrak{e}\in[0,1)$ and $\beta>0$, the point $\alpha>0$ satisfying $\nu_{\omega}(\xi_{\alpha,\beta,\mathfrak{e}})\neq0$ is isolated.
(3) For any fixed $\mathfrak{e}\in[0,1)$ and $k\in\mathbb{R}$, the point $\alpha>0$ satisfying $\nu_{\omega}(\xi_{\alpha,k\alpha,\mathfrak{e}})\neq0$ is isolated.
\end{cor}
\begin{proof}
We just proof $(1)$, the proof of $(2), (3)$ are similar. For any fixed $\mathfrak{e}\in[0,1)$ and $\alpha>0$, if the point $\beta>0$ satisfying $\nu_{\omega}((\xi_{\alpha,\beta,\mathfrak{e}})\neq0$ is not isolated, then there exist an nonzero sequence $\epsilon_{l}\to0$ such that $\beta+\epsilon_{l}>0$ and $\nu_{\omega}(\xi_{\alpha,\beta+\epsilon_{l},\mathfrak{e}})\neq0$. Without loss of generality, suppose that $\{\epsilon_{l}\}$ is a strictly decreasing sequence and $\epsilon_{l}\to0^{+}$, then
$$
\begin{aligned}
&i_{\omega}(\xi_{\alpha,\beta+\epsilon_{1}},\mathfrak{e})-i_{\omega}(\xi_{\alpha,\beta},\mathfrak{e})\\
=&\sum_{j=1}^{l-1}
\left(i_{\omega}(\xi_{\alpha,\beta+\epsilon_{j}},\mathfrak{e}))-i_{\omega}(\xi_{\alpha,\beta+\epsilon_{j+1}},\mathfrak{e})\right)
+i_{\omega}(\xi_{\alpha,\beta+\epsilon_{l}},\mathfrak{e})-i_{\omega}(\xi_{\alpha,\beta},\mathfrak{e}).
\end{aligned}
$$
From Proposition \ref{monot}, for any $j=1,...l-1,$ we have
$$
i_{\omega}(\xi_{\alpha,\beta+\epsilon_{j}},\mathfrak{e})-i_{\omega}(\xi_{\alpha,\beta+\epsilon_{j+1}},\mathfrak{e})
\geq\nu_{\omega}(\xi_{\alpha,\beta+\epsilon_{j+1}},\mathfrak{e})),
$$
hence
$$
i_{\omega}(\xi_{\alpha,\beta+\epsilon_{1}},\mathfrak{e})-i_{\omega}(\xi_{\alpha,\beta},\mathfrak{e})
\geq\sum_{j=1}^{l-1}\nu_{\omega}(\xi_{\alpha,\beta+\epsilon_{j+1}},\mathfrak{e}),
$$
let $l\to+\infty$, then $i_{\omega}(\xi_{\alpha,\beta+\epsilon_{1}},\mathfrak{e})=+\infty$, which contradict with $i_{\omega}(\xi_{\alpha,\beta+\epsilon_{1}},\mathfrak{e})<+\infty$,
hence the point $\beta>0$ satisfying $\nu(\xi_{\alpha,\beta,\mathfrak{e}})\neq0$ is isolated.
\end{proof}
In the following, we give some lemmas which study the Maslov index of $\xi_{\alpha,\beta,\mathfrak{e}}$ for $(\alpha,\beta)$ in some special region, it's important in our
proof of the main theorem.
\begin{lem}\cite[Proposition 3.6 and Remark 3.7 ]{HS}\label{Lag} When $\alpha=1/2,\beta=3/2$, it is just the linear system corresponding the Kepler orbit, and
for any $\omega\in\mathbb{U}\setminus\{1\}$, $\mathfrak{e}\in[0,1)$,
\bea
&&i_{1}(\xi_{1/2,3/2,\mathfrak{e}})=0,\;\nu_{1}(\xi_{1/2,3/2,\mathfrak{e}})=3,\;\nonumber\\
&&i_{\omega}(\xi_{1/2,3/2,\mathfrak{e}})=2,\;\nu_{\omega}(\xi_{1/2,3/2,\mathfrak{e}})=0.
\eea
\end{lem}
\begin{lem}\label{cri-index2}
When $0<\alpha<1/2, \beta=3\alpha$, we have $i_{1}(\xi_{\alpha,\beta,\mathfrak{e}})=0$ and $v_{1}(\xi_{\alpha,\beta,\mathfrak{e}})=0$.
\end{lem}
\begin{proof}
From Proposition \ref{monot2}, we have $i_{1}(\xi_{1/2,3/2,\mathfrak{e}})\geq i_{1}(\xi_{\alpha,\beta,\mathfrak{e}})+v_{1}(\xi_{\alpha,\beta,\mathfrak{e}})$,
Lemma \ref{Lag} tells us that $i_{1}(\xi_{1/2,3/2,\mathfrak{e}})=0$, hence $i_{1}(\xi_{\alpha,\beta,\mathfrak{e}})=0, v_{1}(\xi_{\alpha,\beta,\mathfrak{e}})=0$.
\end{proof}
\begin{lem}\label{cri-eig}
When $\beta=\alpha+1$, the spectrum of $\xi_{\alpha,\beta,\mathfrak{e}}(2\pi)$ is $\sigma(\xi_{\alpha,\beta,\mathfrak{e}}(2\pi))=\{1\}\cup\sigma(\Phi_{\alpha,\mathfrak{e}}(2\pi))$, where
$\Phi_{\alpha,\mathfrak{e}}(\theta)$ satisfies the following system,
\bea\label{Ham2}
\Phi'_{\alpha,\mathfrak{e}}(\theta)=J_{2}\left(
                     \begin{array}{cc}
                       1 & 0 \\
                       0 & 4-\frac{2+2\alpha}{1+\mathfrak{e}\cos\theta} \\
                     \end{array}
                   \right)\Phi_{\alpha,\mathfrak{e}}(\theta),\;\Phi_{\alpha,\mathfrak{e}}(0)=I_{2}.
\eea
\end{lem}
\begin{proof}
 When $\beta=\alpha+1$, the linear Hamiltonian system (\ref{a-b-type}) becomes
\bea
\xi'_{\alpha,\beta,\mathfrak{e}}(\theta)=J_{4}\mathfrak{B}_{\alpha,\beta,\mathfrak{e}}(\theta)\xi_{\alpha,\beta,\mathfrak{e}}(\theta),\; \xi_{\alpha,\beta,\mathfrak{e}}(0)=I_{4}.\label{equ b=a+1}
\eea
with
\bea \mathfrak{B}_{\alpha,\beta,\mathfrak{e}}(\theta)=\left( \begin{array}{cccc} I_{2} & -J_{2} \\
J_{2} & I_{2}-\left(\begin{array}{cc}\frac{2+2\alpha}{1+\mathfrak{e}\cos\theta} & 0 \\ 0 & 0\end{array}\right)
\end{array}\right).\nonumber
\eea
Let $z(\theta)=(x_{1},x_{2},y_{1},y_{2})^{\mathcal{T}}(\theta)$ satisfies $z^{'}(\theta)=J_{4}\mathfrak{B}_{\alpha,\beta,\mathfrak{e}}(\theta)z(\theta)$, then direct computations shows that
$$
x_{2}(\theta)+y_{1}(\theta)=C_{1},\ \ y_{2}^{'}=-2y_{1}(\theta)+C_{1},\ \
y_{1}^{'}=x_{1}+y_{2};
$$
$$
y_{1}^{''}+(4-\frac{2+2\alpha}{1+\mathfrak{e}\cos\theta})y_{1}(\theta)=2C_{1};
$$
where $C_{1}$ is a constant. One can see that the last equation is just a Hill's equation, and if we know the solutions of this Hill's equation, base on the above relation, we can obtain all the solutions of (\ref{equ b=a+1}). More precisely, we have
\bea\label{fund solution}
\begin{aligned}
&x_{1}(\theta)=y^{'}_{1}(\theta)-C_{1}\theta+2\int_{0}^{\theta}y_{1}(s)ds-C_{2},\ \ x_{2}(\theta)=C_{1}-y_{1}(\theta),\\
&y_{1}(\theta)=y_{1}(\theta),\ \ y_{2}(\theta)=C_{1}\theta-2\int_{0}^{\theta}y_{1}(s)ds+C_{2}.
\end{aligned}
\eea
where $C_{2}$ is a constant. Let $v_{1}(\theta), v_{2}(\theta)$ satisfying $(v_{1}(0), v_{2}(0))=(1,0), (v'_{1}(0), v'_{2}(0))=(0,1)$ be two nonzero linearly independent solutions of the following homogeneous Hill equation,
\bea\label{Hill equ}
y_{1}^{''}+(4-\frac{2+2\alpha}{1+\mathfrak{e}\cos\theta})y_{1}(\theta)=0.
\eea
In (\ref{fund solution}), take $(C_1,C_2,y_1(\theta))$ successively as
$$
 (0,0,v_2(\theta)),\quad (1,0,v_0(\theta)),\quad
 (0,0,v_1(\theta)),\quad (0,1,0),
$$
where $v_0$ is the solution of
$$
 v_0''+\left(4-\frac{2+2\alpha}{1+\mathfrak{e}\cos\theta}\right)v_0=2,
 \qquad v_0(0)=v_0'(0)=0.
$$
By variation of constants,
$$
 v_0(\theta)=2\int_0^\theta
 \bigl(v_2(\theta)v_1(s)-v_1(\theta)v_2(s)\bigr)\,ds.
$$
Substituting these four choices into (\ref{fund solution}) gives
$$
 Z(\theta)=
 \begin{pmatrix}
 v_2'+2\int_0^\theta v_2 & v_0'-\theta+2\int_0^\theta v_0
                           & v_1'+2\int_0^\theta v_1 & -1\\
 -v_2 & 1-v_0 & -v_1 & 0\\
 v_2  & v_0   & v_1  & 0\\
 -2\int_0^\theta v_2 & \theta-2\int_0^\theta v_0
                      & -2\int_0^\theta v_1 & 1
 \end{pmatrix},
$$
where $\int_0^\theta f$ denotes $\int_0^\theta f(s)\,ds$.
Since $\det Z(0)=1$, the four columns are linearly independent and
the fundamental solution of (\ref{equ b=a+1}) is
$$
\begin{aligned}
 &\xi_{\alpha,\beta,\mathfrak{e}}(\theta)=Z(\theta)Z(0)^{-1}\\
 &\quad=
 \begin{pmatrix}
 v_2'+2\int_0^\theta v_2 &
 v_0'-\theta+2\int_0^\theta v_0 &
 v_0'+v_1'-\theta+2\int_0^\theta (v_0+v_1) &
 v_2'-1+2\int_0^\theta v_2\\
 -v_2 & 1-v_0 & 1-v_0-v_1 & -v_2\\
 v_2 & v_0 & v_0+v_1 & v_2\\
 -2\int_0^\theta v_2 &
 \theta-2\int_0^\theta v_0 &
 \theta-2\int_0^\theta (v_0+v_1) &
 1-2\int_0^\theta v_2
 \end{pmatrix}.
\end{aligned}
$$
To determine its spectrum, consider
$$
 P\xi_{\alpha,\beta,\mathfrak{e}}(\theta)P^{-1}
 =\begin{pmatrix}
 1 & \theta-2\int_0^\theta v_0 & -2\int_0^\theta v_2
                                & -2\int_0^\theta v_1\\
 0 & 1 & 0 & 0\\
 0 & v_0' & v_2' & v_1'\\
 0 & v_0  & v_2  & v_1
 \end{pmatrix},\ \ P=\begin{pmatrix}
 0&0&0&1\\0&1&1&0\\1&0&0&1\\0&0&1&0
 \end{pmatrix},
$$
then we have
$$
 \det\!\left(\lambda I_4-\xi_{\alpha,\beta,\mathfrak{e}}(\theta)\right)
 =(\lambda-1)^2
 \det\!\left(\lambda I_2-\Phi_{\alpha,\mathfrak{e}}(\theta)\right).
$$
Therefore, apart from the eigenvalue $1$, the remaining eigenvalues of the fundamental solution matrix $\xi_{\alpha,\beta,\mathfrak{e}}(\theta)$
are determined by the $2\times2$ fundamental solution $\Phi_{\alpha,\mathfrak{e}}(\theta)$ which satisfies
$$
\Phi_{\alpha,\mathfrak{e}}(\theta)=\left(
                                                                      \begin{array}{cc}
                                                                        v'_{2}(\theta) & v'_{1}(\theta) \\
                                                                        v_{2}(\theta) & v_{1}(\theta) \\
                                                                      \end{array}
                                                                    \right),\ \
\Phi'_{\alpha,\mathfrak{e}}(\theta)=J_{2}\left(
                     \begin{array}{cc}
                       1 & 0 \\
                       0 & 4-\frac{2+2\alpha}{1+\mathfrak{e}\cos\theta} \\
                     \end{array}
                   \right)\Phi_{\alpha,\mathfrak{e}}(\theta),\;\Phi_{\alpha,\mathfrak{e}}(0)=I_{2},
$$
hence $\sigma(\xi_{\alpha,\beta,\mathfrak{e}}(2\pi))=\{1\}\cup\sigma(\Phi_{\alpha,\mathfrak{e}}(2\pi))$.
\end{proof}
\begin{lem}\label{cri-index}
When $\beta=\alpha+1$, we have
\begin{description}
  \item[(1)] If $\alpha=\frac{1}{2}$, then for any $\mathfrak{e}\in[0,1)$, we have $i_{1}(\xi_{\alpha,\beta,\mathfrak{e}})=0,\;\nu_{1}(\xi_{\alpha,\beta,\mathfrak{e}})=3$;
  \item[(2)] If $\alpha>\frac{1}{2}$, then for any $\mathfrak{e}\in[0,1)$, we have $i_{1}(\xi_{\alpha,\beta,\mathfrak{e}})=0,\;\nu_{1}(\xi_{\alpha,\beta,\mathfrak{e}})=1$;
  \item[(3)] If $-1<\alpha<\frac{1}{2}$, then for any $\mathfrak{e}\in[0,1)$, we have $i_{1}(\xi_{\alpha,\beta,\mathfrak{e}})=2,\;\nu_{1}(\xi_{\alpha,\beta,\mathfrak{e}})=1$.
  \item[(4)] For any $\mathfrak{e}\in[0,1)$, $\omega\in\mathbb{U}$, if $\alpha>1+\frac{3}{2}\mathfrak{e}$, then $i_{\omega}(\xi_{\alpha,\beta,\mathfrak{e}})=0$ and $\xi_{\alpha,\beta,\mathfrak{e}}(2\pi)\approx N_{1}(1,-1)\diamond D(\lambda)$ with some $1\neq\lambda>0$.
\end{description}
\end{lem}
\begin{rem}\label{rem (4)}
By the semi-continuity of the Morse index, \textbf{(4)} of Lemma \ref{cri-index} also implies $i_{\omega}(\xi_{\alpha_{0},\alpha_{0}+1,\mathfrak{e}})=0$ where $\alpha_{0}=1+3\mathfrak{e}/2$.
\end{rem}
\begin{proof}
(1)\;For $\alpha=1/2$, this just follows from Lemma \ref{Lag}. In fact, in the recent paper \cite{OS}, Ou and Sun further obtain
$\ker(\mathcal{A}(\mathfrak{R}_{1/2,3/2},\mathfrak{e}))=\text{Span}_{\mathbb{C}}\{z_{1}(\theta),z_{2}(\theta),z_{3}(\theta)\}$, where
\bea
z_{1}(\theta)&=&(\sin{\theta}+\mathfrak{e}\sin{\theta}\cos{\theta}\ , \ 2\cos{\theta}+\mathfrak{e}\cos^{2}{\theta})^\mathcal{T},\label{solu1} \\
z_{2}(\theta)&=&(0\ , \ 1)^\mathcal{T},\label{solu2}\\
z_{3}(\theta)&=&(\cos{\theta}+\mathfrak{e}\cos^{2}{\theta}\ , \ -2\sin{\theta}-\mathfrak{e}\sin{\theta}\cos{\theta})^\mathcal{T}.\label{solu3}
\eea

(2)\; For $\alpha>\frac{1}{2}$, we have
$\mathcal{A}(\mathfrak{R}_{\alpha,\beta},\mathfrak{e})\geq\mathcal{A}(\mathfrak{R}_{1/2,3/2},\mathfrak{e})$, this combine with Proposition \ref{rela-Morse-Maslov} and case $(1)$, we know
$\phi_{1}(\mathcal{A}(\mathfrak{R}_{1/2,3/2},\mathfrak{e}))=i_{1}(\xi_{1/2,3/2,\mathfrak{e}})=0,$ hence $$\phi_{1}(\mathcal{A}(\mathfrak{R}_{\alpha,\beta},\mathfrak{e}))=i_{1}(\xi_{\alpha,\beta,\mathfrak{e}})=0.$$ On the other hand, if $0\neq x=(x_{1},x_{2})^{\mathcal{T}}$ satisfies
$\mathcal{A}(\mathfrak{R}_{\alpha,\beta},\mathfrak{e})x=0$, then it's equivalent to
$$
\mathcal{A}(\mathfrak{R}_{1/2,3/2},\mathfrak{e})x=-
\frac{2\alpha-1}{1+\mathfrak{e}\cos\theta}\left(\begin{array}{cc}
1 & 0 \\
0 & 0 \\
\end{array}
\right)x.
$$
Since $\mathcal{A}(\mathfrak{R}_{1/2,3/2},\mathfrak{e})\geq0$, the only possible case is $x=(0,x_{2})^{\mathcal{T}}$ and $\mathcal{A}(\mathfrak{R}_{1/2,3/2},\mathfrak{e})x=0$, but direct computation shows that
$$
0=\left(\mathcal{A}(\mathfrak{R}_{1/2,3/2},\mathfrak{e})x,x\right)_{L^{2}}
=\int_{0}^{2\pi}|\dot{x}_{2}(t)|^{2}dt.
$$
The last equality is obtained by integration by parts. Therefore, $x_{2}(t)$ must be a non-zero constant, then $$\nu_{1}(\mathcal{A}(\mathfrak{R}_{\alpha,\beta},\mathfrak{e}))=
\nu_{1}(\xi_{\alpha,\beta,\mathfrak{e}})=1.$$

(3)\; For $-1<\alpha<\frac{1}{2}$, consider $z_{1}(\theta),z_{2}(\theta)$ and $z_{3}(\theta)$ given by (\ref{solu1}), (\ref{solu2}) and (\ref{solu3}), then from
$$
\mathcal{A}(\mathfrak{R}_{\alpha,\beta},\mathfrak{e})=\mathcal{A}(\mathfrak{R}_{1/2,3/2},\mathfrak{e})+
\frac{2\alpha-1}{1+\mathfrak{e}\cos\theta}\left(\begin{array}{cc}
1 & 0 \\
0 & 0 \\
\end{array}
\right),
$$
we have
\bea\label{equ1}
\mathcal{A}(\mathfrak{R}_{\alpha,\beta},\mathfrak{e})z_{2}=0,\;\text{and}\;\left(\mathcal{A}(\mathfrak{R}_{\alpha,\beta},\mathfrak{e})z_{i},z_{i}\right)_{L^{2}}<0\;\text{for}\; i=1,3,
\eea
hence $\phi_{1}(\mathcal{A}(\mathfrak{R}_{\alpha,\beta},\mathfrak{e}))\geq2$. On the other hand,
$$
R(\theta)\mathcal{A}(\mathfrak{R}_{-1,0},\mathfrak{e})R^{\mathcal{T}}(\theta)=R(\theta)(-\frac{d^2}{d\theta^2}I_{2}-2J_{2}\frac{d}{d\theta})R^{\mathcal{T}}(\theta)
=-\frac{d^2}{d\theta^2}I_{2}-I_{2},
$$
it's Morse index $\phi_{1}(\mathcal{A}(\mathfrak{R}_{-1,0},\mathfrak{e}))=2$, this combine with the fact $\mathcal{A}(\mathfrak{R}_{\alpha,\beta},\mathfrak{e})\geq
\mathcal{A}(\mathfrak{R}_{-1,0},\mathfrak{e})$ for $-1<\alpha<1/2$, which implies $$\phi_{1}(\mathcal{A}(\mathfrak{R}_{\alpha,\beta},\mathfrak{e}))\leq2.$$ Therefore, we finally obtain
\bea\label{equ2}
\phi_{1}(\mathcal{A}(\mathfrak{R}_{\alpha,\beta},\mathfrak{e}))=i_{1}(\xi_{\alpha,\beta,\mathfrak{e}})=2,\;\text{for any}\; -1<\alpha<1/2
\eea
Last, we show that $\nu(\xi_{\alpha,\beta,\mathfrak{e}})=1$. Similar, we consider $0\neq x=(x_{1},x_{2})^{\mathcal{T}}$ which satisfies
$\mathcal{A}(\mathfrak{R}_{\alpha,\beta},\mathfrak{e})x=0$, then take $\tilde{\alpha}\in(-1,\alpha)$, we have
$$
\mathcal{A}(\mathfrak{R}_{\tilde{\alpha},\tilde{\alpha}+1},\mathfrak{e})=\mathcal{A}(\mathfrak{R}_{\alpha,\beta},\mathfrak{e})-
\frac{2(\alpha-\tilde{\alpha})}{1+\mathfrak{e}\cos\theta}\left(\begin{array}{cc}
1 & 0 \\
0 & 0 \\
\end{array}
\right).
$$
If $x_{1}\neq0$, then $<\mathcal{A}(\mathfrak{R}_{\tilde{\alpha},\tilde{\alpha}+1},\mathfrak{e})x,x>_{L^{2}}<0$. Also, (\ref{equ1}) implies
$<\mathcal{A}(\mathfrak{R}_{\tilde{\alpha},\tilde{\alpha}+1},\mathfrak{e})z_{i},z_{i}>_{L^{2}}<0$ for $i=1,3,$ hence $
\phi_{1}(\mathcal{A}(\mathfrak{R}_{\tilde{\alpha},\tilde{\alpha}+1},\mathfrak{e}))\geq3,
$
which contradict with $(\ref{equ2})$. This means $x_{1}=0$, and $x=(0,x_{2})^{\mathcal{T}}$, then similar the proof of $(2)$,
$\mathcal{A}(\mathfrak{R}_{\alpha,\beta},\mathfrak{e})x=0$ is equivalent to
$$
\mathcal{A}(\mathfrak{R}_{1/2,3/2},\mathfrak{e})x=-
\frac{2\alpha-1}{1+\mathfrak{e}\cos\theta}\left(\begin{array}{cc}
1 & 0 \\
0 & 0 \\
\end{array}
\right)x,
$$
then
$$
0=\left(\mathcal{A}(\mathfrak{R}_{1/2,3/2},\mathfrak{e})x,x\right)_{L^{2}}
=\int_{0}^{2\pi}|\dot{x}_{2}(t)|^{2}dt,
$$
which implies $x_{2}(t)$ must be a non-zero constant, hence
$$\nu_{1}(\mathcal{A}(\mathfrak{R}_{\alpha,\beta},\mathfrak{e}))=
\nu_{1}(\xi_{\alpha,\beta,\mathfrak{e}})=1.$$

(4) From Lemma \ref{cri-eig}, we know $\sigma(\xi_{\alpha,\beta,\mathfrak{e}}(2\pi))=\{1\}\cup\sigma(\Phi_{\alpha,\mathfrak{e}}(2\pi))$, where $\Phi_{\alpha,\mathfrak{e}}(\theta)$ satisfies the linear equation (\ref{Ham2}). This equation corresponding to Hill's equation (\ref{Hill equ}) and the corresponding Sturm-Liouville operator is
$$
\mathcal{F}_{\alpha,\mathfrak{e}}=-\frac{d^2}{d\theta^2}-4+\frac{2+2\alpha}{1+\mathfrak{e}\cos\theta}\;,
$$
which is a self-adjoint operator in $L^2([0,2\pi], \mathbb{C})$ with domain
$$ \bar{D}_{1}(\omega,2\pi)=\{y\in W^{2,2}([0,2\pi],\mathbb{C})\;|\;y(2\pi)=\omega y(0),\dot{y}(2\pi)=\omega \dot{y}(0)\},\;\omega\in\mathbb{U}. $$

For $\alpha>1+\frac{3}{2}\mathfrak{e}$ and any $\tilde{\mathfrak{e}}\in[0,\mathfrak{e}]$, we have $\frac{1+2\alpha}{1+\tilde{\mathfrak{e}}\cos\theta}>3$, hence
$
\mathcal{F}_{\alpha,\tilde{\mathfrak{e}}}>A_{\tilde{\mathfrak{e}}},
$
where $A_{\tilde{\mathfrak{e}}}=-\frac{d^2}{d\theta^2}-1+\frac{1}{1+\tilde{\mathfrak{e}}\cos\theta}$. Since $R(\theta)\mathcal{A}(I_{2},\tilde{\mathfrak{e}})R^{\mathcal{T}}(\theta)=A_{\tilde{\mathfrak{e}}}\oplus A_{\tilde{\mathfrak{e}}}$, Lemma \ref{lem-pos} and Lemma \ref{lem-pos2} imply $A_{\tilde{\mathfrak{e}}}\geq0$ in $\bar{D}_{1}(\omega,2\pi)$ for any $\omega\in\mathbb{U}$, hence $\mathcal{F}_{\alpha,\tilde{\mathfrak{e}}}>0$
in $\bar{D}_{1}(\omega,2\pi)$. Similar to Proposition \ref{rela-Morse-Maslov}, we also have the relation $i_{\omega}(\Phi_{\alpha,\tilde{\mathfrak{e}}})=\phi_{\omega}(\mathcal{F}_{\alpha,\tilde{\mathfrak{e}}})=0,$ and $\nu_{\omega}(\Phi_{\alpha,\tilde{\mathfrak{e}}})=\nu_{\omega}(\mathcal{F}_{\alpha,\tilde{\mathfrak{e}}})=0$, this implies any $\omega\in\mathbb{U}$ is not the eigenvalue of $\Phi_{\alpha,\tilde{\mathfrak{e}}}(2\pi)$, hence $\Phi_{\alpha,\tilde{\mathfrak{e}}}(2\pi)$ is always hyperbolic for any $\tilde{\mathfrak{e}}\in[0,\mathfrak{e}]$. Since $\Phi_{\alpha,0}(2\pi)$ has eigenvalues $e^{\pm2\pi\sqrt{2\alpha-2}}$, we have $\Phi_{\alpha,\tilde{\mathfrak{e}}}(2\pi)\approx D(\lambda)$ with some $1\neq\lambda>0$ for any $\tilde{\mathfrak{e}}\in[0,\mathfrak{e]}$, then $\xi_{\alpha,\alpha+1,\tilde{\mathfrak{e}}}(2\pi)\approx N_{1}(1,a)\diamond D(\lambda)$ with $1\neq\lambda>0$.
Moreover, when $\alpha>1+\frac{3}{2}\mathfrak{e}$, we obtain
\bea\label{inequ1}
\mathcal{A}(\mathfrak{R}_{\alpha,\alpha+1},0)\geq\mathcal{A}(\mathfrak{R}_{7/8,15/8
},0)=\mathcal{A}(\mathfrak{R}_{1/2,3/2},0)+
\frac{3}{4}\left(\begin{array}{cc}
1 & 0 \\
0 & 0 \\
\end{array}
\right),
\eea
and $\mathcal{A}(\mathfrak{R}_{7/8,15/8},0)$ in $D_{2}(-1,2\pi)$ has
$$
\ker\mathcal{A}(\mathfrak{R}_{7/8,15/8},0)=\text{Span}_{\mathbb{C}}\left\{\left(\begin{array}{c}
                                                                     \cos\frac{\theta}{2} \\
                                                                     -4\sin\frac{\theta}{2}
                                                                   \end{array}\right),\;\left(\begin{array}{c}
                                                                     \sin\frac{\theta}{2} \\
                                                                     4\cos\frac{\theta}{2}
                                                                   \end{array}\right)\right\},
$$
this combine with (\ref{inequ1}), we obtain
$$
\phi_{-1}(\mathcal{A}(\mathfrak{R}_{1/2,3/2},0))\geq\phi_{-1}(\mathcal{A}(\mathfrak{R}_{7/8,15/8},0))+\nu_{-1}(\mathcal{A}(\mathfrak{R}_{7/8,15/8},0))
=\phi_{-1}(\mathcal{A}(\mathfrak{R}_{7/8,15/8},0))+2.
$$
On the other hand, Lemma \ref{Lag} tells us that $\phi_{-1}(\mathcal{A}(\mathfrak{R}_{1/2,3/2},0))=i_{-1}(\xi_{1/2,3/2,0})=2$, then combine with above inequality, we obtain $\phi_{-1}(\mathcal{A}(\mathfrak{R}_{7/8,15/8,0})=0,$
and hence
$
i_{-1}(\xi_{\alpha,\alpha+1,0})=\phi_{-1}(\mathcal{A}(\mathfrak{R}_{\alpha,\alpha+1},0))=0.
$
From \textbf{(2)}, we know $i_{1}(\xi_{\alpha,\alpha+1,0})=0$, this combine with the Bott iteration formula in Proposition \ref{Bott} and the Proposition \ref{split norm} of
the Splitting number, we have
\bea\label{equal2}
0=i_{-1}(\xi_{\alpha,\alpha+1,0})=i_{1}(\xi_{\alpha,\alpha+1,0})+S_{M}^{+}(1)+
\sum_{0<\theta<\vartheta}(S_{M}^{+}(e^{\sqrt{-1}\theta})-S_{M}^{-}(e^{\sqrt{-1}\theta}))\nonumber
-S_{M}^{-}(-1)=S_{M}^{+}(1),
\eea
where $M=\xi_{\alpha,\alpha+1,0}(2\pi)$, the only possible case is $\xi_{\alpha,\alpha+1,0}(2\pi)\approx N_{1}(1,-1)\diamond D(\lambda)$ with some $1\neq\lambda>0$. Since we always have
$\nu_{1}(\xi_{\alpha,\alpha+1,\tilde{\mathfrak{e}}})=1$, we obtain $\xi_{\alpha,\alpha+1,\tilde{\mathfrak{e}}}(2\pi)\approx N_{1}(1,-1)\diamond D(\lambda)$ for any $\tilde{\mathfrak{e}}\in[0,\mathfrak{e}]$. This combine with the fact $i_{1}(\xi_{\alpha,\beta,\tilde{\mathfrak{e}}})=0$ in \textbf{(2)}, using the Bott iteration formula
again, we have $i_{\omega}(\xi_{\alpha,\alpha+1,\tilde{\mathfrak{e}}})=0$ for any $\omega\in\mathbb{U}$ and $\tilde{\mathfrak{e}}\in[0,\mathfrak{e}]$.
This completes the proof.
\end{proof}
\begin{lem}\label{even-nullity}
When $\beta\neq\alpha+1$, $\nu_{1}(\xi_{\alpha,\beta,\tilde{\mathfrak{e}}})$ must be even number.
\end{lem}
\begin{proof}
This essentially follows from Lemma $8$ in \cite{MSS1}, which says that if $0\neq x=(x_{1},x_{2})^{\mathcal{T}}\in\text{Ker}(\mathcal{A}(\mathfrak{R}_{\alpha,\beta},\mathfrak{e}))$, where
$$
\begin{aligned}
&x_{1}(\theta)=a_{0}+\sum_{k\geq1}a_{k}\cos(k\theta)+\sum_{k\geq1}b_{k}\sin(k\theta),\\
&x_{2}(\theta)=c_{0}+\sum_{k\geq1}c_{k}\cos(k\theta)+\sum_{k\geq1}d_{k}\sin(k\theta),
\end{aligned}
$$
then one can construct two independent solutions $x^{(1)}=(x^{(1)}_{1},x^{(1)}_{2})^{\mathcal{T}}, x^{(2)}=(x^{(2)}_{1},x^{(2)}_{2})^{\mathcal{T}}\in\text{Ker}(\mathcal{A}(\mathfrak{R}_{\alpha,\beta},\mathfrak{e}))$ with
$$
\begin{aligned}
&x_{1}^{(1)}(\theta)=a_{0}+\sum_{k\geq1}a_{k}\cos(k\theta),\; x_{2}^{(1)}(\theta)=\sum_{k\geq1}d_{k}\sin(k\theta),\\
&x_{1}^{(2)}(\theta)=\sum_{k\geq1}a_{k}\sin(k\theta),\;x_{2}^{(2)}(\theta)=c_{0}-\sum_{k\geq1}d_{k}\cos(k\theta),
\end{aligned}
$$
where $a_{0}=-\mathfrak{e}(1+\alpha+\beta)^{-1}(a_{1}/2+d_{1})$, $c_{0}=\mathfrak{e}(1+\alpha-\beta)^{-1}(a_{1}+d_{1}/2)$. Moreover, one can check that
for any another solution $\tilde{x}=(\tilde{x}_{1},\tilde{x}_{2})^{\mathcal{T}}\in\text{Ker}(\mathcal{A}(\mathfrak{R}_{\alpha,\beta},\mathfrak{e}))$, let $\tilde{x}^{(1)}$ and
 $\tilde{x}^{(2)}$ denote the linearly independent solutions constructed corresponding to $\tilde{x}$, then either $\text{Span}\{x^{(1)},x^{(2)}\}=\text{Span}\{\tilde{x}^{(1)},\tilde{x}^{(2)}\}$, or $\text{Span}\{x^{(1)},x^{(2)}\}\cap\text{Span}\{\tilde{x}^{(1)},\tilde{x}^{(2)}\}=\{0\}$.
 Therefore, $\text{Ker}(\mathcal{A}(\mathfrak{R}_{\alpha,\beta},\mathfrak{e}))$ has even dimension, then $\nu_{1}(\xi_{\alpha,\beta,\mathfrak{e}})$ must be even number.
\end{proof}
\begin{cor}\label{identity-block}
Let $\mathfrak{e}\in[0,1)$, $\alpha>0$ and $\beta>\alpha+1$. If $\nu_{1}(\xi_{\alpha,\beta,\mathfrak{e}})\neq0$, then
$$
\xi_{\alpha,\beta,\mathfrak{e}}(2\pi)\approx I_{2}\diamond M_{2}
\quad\text{for some }M_{2}\in\Sp(2).
$$
\end{cor}
\begin{proof}
Since $\nu_{1}(\xi_{\alpha,\beta,\mathfrak{e}})>0$, let $x=(x_{1},x_{2})^{\mathcal{T}}\in\text{Ker}(\mathcal{A}(\mathfrak{R}_{\alpha,\beta},\mathfrak{e}))$, and
take $x^{(1)},x^{(2)}$ be the two linearly independent periodic
solutions constructed in the proof of Lemma \ref{even-nullity}.
The corresponding solutions of the first-order system are
$$
z_{j}=\binom{(x^{(j)})'+J_{2}x^{(j)}}{x^{(j)}},
\qquad j=1,2.
$$
Take $X^{(j)}(\theta)=(1+\mathfrak{e}\cos\theta)^{-1}x^{(j)}(\theta)$,
direct computation shows that
\bea\label{equa15}
(1+\mathfrak{e}\cos\theta)^{-2}z_{1}^{\mathcal{T}}J_{4}z_{2}
=\Bigl((X^{(1)})^{\mathcal{T}}(X^{(2)})'
-((X^{(1)})')^{\mathcal{T}}X^{(2)}
+2(X^{(1)})^{\mathcal{T}}J_{2}X^{(2)}\Bigr).
\eea
Moreover, the equation
$(x^{(j)})''+2J_{2}(x^{(j)})'=\mathfrak{R}_{\alpha,\beta}X^{(j)}
$ gives
$$
X^{(1)}=\sum_{k\geq1}\binom{A_{k}\cos(k\theta)}{B_{k}\sin(k\theta)},
\qquad
X^{(2)}=\sum_{k\geq1}\binom{A_{k}\sin(k\theta)}{-B_{k}\cos(k\theta)},
$$
where 
\bea\label{Fourier}
(1+\alpha+\beta)A_{k}=-k^{2}a_{k}-2kd_{k},\ \
(1+\alpha-\beta)B_{k}=-k^{2}d_{k}-2ka_{k}.
\eea. Substituting this Fourier expansion into equation \ref{equa15} and integrating both sides over a period of
$2\pi$, we obtain
$$
\begin{aligned}
\frac{1}{2\pi}\int_{0}^{2\pi}(1+\mathfrak{e}\cos\theta)^{-2}z_{1}^{\mathcal{T}}J_{4}z_{2}d\theta&=\sum_{k\geq1}\bigl[k(A_{k}^{2}+B_{k}^{2})+2A_{k}B_{k}\bigr]\\
 &=\sum_{k\geq1}\bigl[(k-1)(A_{k}^{2}+B_{k}^{2})
                       +(A_{k}+B_{k})^{2}\bigr]\geq0.
\end{aligned}
$$
If the right hand side  equals zero, then $A_{k}=B_{k}=0$ for $k\geq2$ and $B_{1}=-A_{1}$. Thus
$$
x^{(1)}(\theta)=(1+\mathfrak{e}\cos\theta)A_{1}(\cos\theta,-\sin\theta)^{\mathcal{T}}.
$$
Recall that $x^{(1)}=(x_{1}^{(1)},x_{2}^{(1)})^{\mathcal{T}}$, where
$$
\begin{aligned}
&x_{1}^{(1)}(\theta)=a_{0}+\sum_{k\geq1}a_{k}\cos(k\theta),\; x_{2}^{(1)}(\theta)=\sum_{k\geq1}d_{k}\sin(k\theta),
%&x_{1}^{(2)}(\theta)=\sum_{k\geq1}a_{k}\sin(k\theta),\;x_{2}^{(2)}(\theta)=c_{0}-\sum_{k\geq1}d_{k}\cos(k\theta),
\end{aligned}
$$
Comparing the Fourier coefficients gives $a_{1}=A_{1}$ and $d_{1}=-A_{1}$.
Consequently, combine with (\ref{Fourier}), we have $(1+\alpha+\beta)A_{1}=-a_{1}-2d_{1}=A_{1}$, then
$A_{1}=0$ and hence $x^{(1)}=0$, a contradiction.
Therefore 
$$0<\int_{0}^{2\pi}(1+\mathfrak{e}\cos\theta)^{-2}z_{1}^{\mathcal{T}}(\theta)J_{4}z_{2}(\theta)d\theta
=\int_{0}^{2\pi}(1+\mathfrak{e}\cos\theta)^{-2}z_{1}^{\mathcal{T}}(0)J_{4}z_{2}(0)d\theta,$$
where the second equality from the face $z_{j}(\theta)=\xi_{\alpha,\beta,\mathfrak{e}}(\theta)z_{j}(0)$, this implies that
$z_{1}^{\mathcal{T}}(0)J_{4}z_{2}(0)\neq0$, and hence $z_{1}(0),z_{2}(0)$
span a symplectic plane. Since $z_{1}(0), z_{2}(0)\in\ker(\xi_{\alpha,\beta,\mathfrak{e}}(2\pi)-I_{4})$, then $\text{Span}\{z_{1}(0),z_{2}(0)\}$
be a two-dimensional symplectic invariant subspace of $\xi_{\alpha,\beta,\mathfrak{e}}(2\pi)$. Then there exists a symplectic matrix $P$ such that
$P^{-1}\xi_{\alpha,\beta,\mathfrak{e}}(2\pi)P=I_{2}\diamond M_{2}$ with $M_{2}\in\Sp(2)$.
\end{proof}
\begin{cor}\label{odd-index}
For any $\mathfrak{e}\in[0,1)$, $\alpha>0$ and $\beta>\alpha+1$, $i_{1}(\xi_{\alpha,\beta,\mathfrak{e}})$ is an odd number.
\end{cor}
\begin{proof}
First, Proposition \ref{monot} and Corollary \ref{isolate} implies that
$$
i_{1}(\xi_{\alpha,\beta,\mathfrak{e}})=i_{1}(\xi_{\alpha,\alpha+1,\mathfrak{e}})+\nu_{1}(\xi_{\alpha,\alpha+1,\mathfrak{e}})+\sum_{\tilde{\beta}\in(\alpha+1,\beta)}\nu_{1}(\xi_{\alpha,\tilde{\beta},\mathfrak{e}}).
$$
Second, Lemma \ref{cri-index}  and Lemma \ref{even-nullity} tell us that $i_{1}(\xi_{\alpha,\alpha+1,\mathfrak{e}})$ is even, $\nu_{1}(\xi_{\alpha,\alpha+1,\mathfrak{e}})$
is odd and $\nu_{1}(\xi_{\alpha,\tilde{\beta},\mathfrak{e}})$ is even for $\tilde{\beta}>\alpha+1$, hence $i_{1}(\xi_{\alpha,\beta,\mathfrak{e}})$ is an odd number.
\end{proof}
In order to proof Theorem \ref{stable-thm}, we also need some important lemmas from \cite{HW} concerning linear stability.
%\begin{lem}(See Lemma 4.11 and and Proposition 4.12 in \cite{HS})\label{instable-cri1}
%If $M\in\Sp(2n)$ is linearly stable, then for  sufficiently small $\epsilon>0$, we have
%$$
%\det(e^{-\epsilon J_{n}}M-I_{2n})>0.
%$$
%Moreover, let $\gamma$ be a symplectic path from $I$ to $M$, then
%$\det(e^{-\epsilon J_{n}}M-I_{2n})>0$ if and only if $i_{1}(\gamma)+n$ is even.
%\end{lem}
\begin{lem}(See Theorem 1.3 and Theorem 1.5 in \cite{HW})\label{instable-cri2}
Let $W\in\Sp(2k)$ satisfying $\det(W-I_{2k})\neq0$, consider a continuously differentiable symplectic path $\gamma$ from $I$ to $W$, then
$\det(W-I_{2k})<0$ if and only if $i_{1}(\gamma)+k$ is odd. Moreover, $\det(W-I_{2k})<0$  implies that there exists a positive eigenvalue
of $W$ which is bigger than $1$.
\end{lem}
\begin{proof}
This lemma essentially originates from Theorems 1.3 and Theorem 1.5 in \cite{HW}. For the reader's convenience, we give a brief explanation. For this symplectic path $\gamma(t)$ with $\gamma(0)=I$ and $\gamma(\tau)=W$, one can find a symmetric matrix path $B:[0,\tau]\mapsto\text{Sym}(2n)$ such that $\dot{\gamma}(t)=JB(t)\gamma(t)$.
From Theorem 1.3 and Theorem 1.5 in \cite{HW}, we know that $\det(W-I_{2k})<0$  if and only if the relative Morse index $I(A,A-B)$ is odd, see the definition of the relative
Morse index in Section \ref{App-index-relation}. Also, from the Proposition \ref{index relation}, we know
$I(A,A-B)=i_{1}(\gamma)+k$, hence $\det(W-I_{2k})<0$ if and only if $i_{1}(\gamma)+k$ is odd. Moreover,  since $W\in\Sp(2k)$,
if $\lambda$ is an eigenvalue of $W$, then so are $\lambda^{-1},\bar{\lambda}$ and $\bar{\lambda}^{-1}$ (some of them may be equal
to each other), and they have the same geometric and algebraic multiplicities. If $1$ is not an
eigenvalue of $W$, then $\det(W-I_{2k})<0$  implies that there exists a positive eigenvalue
of $W$ which is bigger than $1$.
\end{proof}
%\begin{lem}(See Lemma 4.6 and Proposition 4.12 in \cite{HS})\label{instable-cri3}
%Let $\gamma$ be a symplectic path from $I$ to $M\in\Sp(2n)$, If $i_{1}(\gamma)+n$ is
%odd, then $M$ is linearly unstable.
%\end{lem}
%\textbf{Proof of Theorem \ref{stable-thm}-(I)}): For $\beta>\alpha+1$, $\xi_{\alpha,\beta,\mathfrak{e}}(2\pi)\approx D(\lambda)\diamond M$ with $M\in\Sp(2)$, $\lambda>0.$
Now, we can proof our main Theorem \ref{stable-thm} and Theorem \ref{stable-thm2}
\begin{proof}[\textbf{Proof of Theorem \ref{stable-thm}-(I)}]
From Corollary \ref{odd-index}, for any $\mathfrak{e}\in[0,1)$, $\alpha>0$ and $\beta>\alpha+1$, $i_{1}(\xi_{\alpha,\beta,\mathfrak{e}})$ is an odd number, since $k=2$ in our case, we
have $i_{1}(\xi_{\alpha,\beta,\mathfrak{e}})+k$ is an odd number.

If $\det(\xi_{\alpha,\beta,\mathfrak{e}}(2\pi)-I_{4})\neq0$, then Lemma \ref{instable-cri2} implies there exists a positive eigenvalue
of $\xi_{\alpha,\beta,\mathfrak{e}}(2\pi)$ which is bigger than $1$, then Proposition \ref{homo} tells us that $\xi_{\alpha,\beta,\mathfrak{e}}(2\pi)\approx D(\lambda)\diamond M_{1}$ with $1\neq\lambda>0$ and some $M_{1}\in\Sp(2)$.

If $\det(\xi_{\alpha,\beta,\mathfrak{e}}(2\pi)-I_{4})=0$, Corollary \ref{identity-block} tell us that we can find a symplectic matrix $P$ such that $P^{-1}\xi_{\alpha,\beta,\mathfrak{e}}(2\pi)P=I_{2}\diamond M_{2}$. Now, we assert that  $M_{2}$ must be $D(\lambda)$ with $1\neq\lambda>0$,
otherwise, $M_{2}$ can only be one of $I_{2}$, $N_{1}(-1,b)$, $R(\vartheta)$ with $\vartheta\in(0,2\pi)\setminus\pi$, or $D(\lambda)$ with $-1\neq\lambda<0$, this implies
$\det(e^{-\epsilon2\pi J_{4}}\xi_{\alpha,\beta,\mathfrak{e}}(2\pi)-I_{4})>0$ for sufficiently small $\epsilon>0$, then Lemma \ref{instable-cri2} implies $i_{1}(\xi^{(\epsilon)}_{\alpha,\beta,\mathfrak{e}})+k$ must be even, where
$\xi^{(\epsilon)}_{\alpha,\beta,\mathfrak{e}}(\theta)=e^{-\epsilon \theta J_{4}  }\xi_{\alpha,\beta,\mathfrak{e}}(\theta)$. From the Definition \ref{def:Maslov-type index} and Remark \ref{homtopy} of the Maslov-type index in Section  \ref{App-index} of Appendix, we know $i_{1}(\xi^{(\epsilon)}_{\alpha,\beta,\mathfrak{e}})=i_{1}(\xi_{\alpha,\beta,\mathfrak{e}})$ for sufficiently small $\epsilon>0$, then $i_{1}(\xi_{\alpha,\beta,\mathfrak{e}})+k$ is also even, contradiction. Therefore, $M_{2}=D(\lambda)$ with some $1\neq\lambda>0$. This completes the proof.
\end{proof}
\begin{rem}
For the case $\nu_{1}(\xi_{\alpha,\beta,\mathfrak{e}}(2\pi))\neq0$, then  Corollary \ref{identity-block} tell us that we can find a symplectic matrix $P$ such that $P^{-1}\xi_{\alpha,\beta,\mathfrak{e}}(2\pi)P=I_{2}\diamond M_{2}$ with $M_{2}\neq N(1,\pm1)$. we can also apply Lemma 3.1 in \cite{ZL} to give another proof showing that $\xi_{\alpha,\beta,\mathfrak{e}}(2\pi)\approx I_{2}\diamond D(\lambda)$ with some $1\neq\lambda>0$. Lemma 3.2 of \cite{ZL} tell us that there exist two path $\xi^{(1)}_{\alpha,\beta,\mathfrak{e}}, \xi^{(2)}_{\alpha,\beta,\mathfrak{e}}\in\mathcal{P}_{2\pi}(2)$ with
$\xi^{(1)}_{\alpha,\beta,\mathfrak{e}}(2\pi)=I_{2}$, $\xi^{(2)}_{\alpha,\beta,\mathfrak{e}}(2\pi)=M_{2}$ such that $\xi_{\alpha,\beta,\mathfrak{e}}$ can homotopy to
$\xi^{(1)}_{\alpha,\beta,\mathfrak{e}}\diamond\xi^{(2)}_{\alpha,\beta,\mathfrak{e}}$ satisfying $i_{1}(\xi_{\alpha,\beta,\mathfrak{e}})=i_{1}(\xi^{(1)}_{\alpha,\beta,\mathfrak{e}})+
i_{1}(\xi^{(2)}_{\alpha,\beta,\mathfrak{e}})$. Since $i_{1}(\xi_{\alpha,\beta,\mathfrak{e}})$ is odd, then one of $i_{1}(\xi^{(1)}_{\alpha,\beta,\mathfrak{e}}),
i_{1}(\xi^{(2)}_{\alpha,\beta,\mathfrak{e}})$ must be even. From the property of Maslov index in Lemma 5.6 of \cite{ZL}, we know $i_{1}(\xi^{(1)}_{\alpha,\beta,\mathfrak{e}})$ is odd, hence $i_{1}(\xi^{(2)}_{\alpha,\beta,\mathfrak{e}})$ is even, the only possible case of $M_{2}$ must be $D(\lambda)$ with some $1\neq\lambda>0$.
\end{rem}
\begin{proof}[\textbf{Proof of Theorem \ref{stable-thm}-(II)}]
%When $\alpha>0$ and $0\leq\beta<\alpha$, $\mathcal{A}(\mathfrak{R}_{\alpha,\beta},\mathfrak{e})$ is strictly positive in $\bar{D}_{2}(\omega,2\pi)$ for any $\omega\in\mathbb{U}$ and
%$\xi_{\alpha,\beta,\mathfrak{e}}(2\pi)$ is hyperbolic.
For any $\mathfrak{e}\in[0,1)$, let $\beta_{0}=\alpha_{0}+1=k(\mathfrak{e})\alpha_{0}$, $\alpha_{0}=1+3\mathfrak{e}/2$, $k(\mathfrak{e})=1+(1+3\mathfrak{e}/2)^{-1}$. First,
consider the case $\alpha\in(0,\alpha_{0}]$, one can see that $k(\mathfrak{e})\alpha=\min\{k(\mathfrak{e})\alpha,\alpha+1\}$. When $\beta<\min\{k(\mathfrak{e})\alpha,\alpha+1\}$, it's equivalent to $\beta<k(\mathfrak{e})\alpha$,
Proposition \ref{monot2} implies
$
i_{\omega}(\xi_{\alpha_{0},\beta_{0},\mathfrak{e}})=i_{\omega}(\xi_{\alpha_{0},k(\mathfrak{e})\alpha_{0},\mathfrak{e}})\geq i_{\omega}(\xi_{\alpha,k(\mathfrak{e})\alpha,\mathfrak{e}})
$.
Moreover, Remark \ref{rem (4)} tells us that $i_{\omega}(\xi_{\alpha_{0},\beta_{0},\mathfrak{e}})=0$, hence $i_{\omega}(\xi_{\alpha,k(\mathfrak{e})\alpha,\mathfrak{e}})=0$ for
any $\omega\in\mathbb{U}$. Using Proposition \ref{monot}, we have
$$
0=i_{\omega}(\xi_{\alpha,k(\mathfrak{e})\alpha,\mathfrak{e}})\geq i_{\omega}(\xi_{\alpha,\beta,\mathfrak{e}})+\nu_{\omega}(\xi_{\alpha,\beta,\mathfrak{e}}),\;\text{for any}\;\beta\in[0,k(\mathfrak{e})\alpha).
$$
This implies that $\nu_{\omega}(\xi_{\alpha,\beta,\mathfrak{e}})=0$ and hence
$\xi_{\alpha,\beta,\mathfrak{e}}(2\pi)$ is hyperbolic for any $\beta\in[0,k(\mathfrak{e})\alpha).$ Second, for the case $\alpha\in[\alpha_{0},+\infty)$, we find
that $\alpha+1=\min\{k(\mathfrak{e}),\alpha+1\}$. Moreover,
\textbf{(4)} of Lemma \ref{cri-index} and Remark \ref{rem (4)} implies $i_{\omega}(\xi_{\alpha,\alpha+1,\mathfrak{e}})=0$.
When $\beta<\min\{k(\mathfrak{e})\alpha,\alpha+1\}$, it's equivalent to $\beta<\alpha+1$,
Proposition \ref{monot} implies
$$
0=i_{\omega}(\xi_{\alpha,\alpha+1,\mathfrak{e}})\geq i_{\omega}(\xi_{\alpha,\beta,\mathfrak{e}})+\nu_{\omega}(\xi_{\alpha,\beta,\mathfrak{e}}),\;\text{for any}\;\beta\in[0,\alpha+1).
$$
This implies that $\nu_{\omega}(\xi_{\alpha,\beta,\mathfrak{e}})=0$ and hence
$\xi_{\alpha,\beta,\mathfrak{e}}(2\pi)$ is hyperbolic for any $\beta\in[0,\alpha+1).$ This completes the proof.
\end{proof}
\begin{proof}[\textbf{Proof of Theorem \ref{stable-thm2}}]
This closely follows the proof of the linear stability of the Lagrangian solution in \cite{HLS}, however, we generalize it to a broader class of systems and a wider range of parameters. First, from Lemma \ref{cri-index2} and Lemma \ref{cri-index}, for any $\mathfrak{e}\in[0,1)$, we have $i_{1}(\xi_{\alpha,\beta,\mathfrak{e}})=0$ if either $\beta=3\alpha$ with $0\leq\alpha\leq1/2$ or $\beta=\alpha+1$ with $\alpha>1/2$. Then from the monotonicity Proposition \ref{monot}, we have $i_{1}(\xi_{\alpha,\beta,\mathfrak{e}})=0, v_{1}(\xi_{\alpha,\beta,\mathfrak{e}})=0$ for all $\mathfrak{e}\in[0,1), \alpha>0, 0\leq\beta<\min\{3\alpha,\alpha+1\}$. Second, from
Lemma \ref{Lag}, we have $i_{-1}(\xi_{1/2,3/2,\mathfrak{e}})=2,\;\nu_{-1}(\xi_{1/2,3/2,\mathfrak{e}})=0$. Therefore, for any $\mathfrak{e}\in[0,1)$, when $\alpha>0, 0\leq\beta<\min\{3\alpha,\alpha+1\}$ and $(\alpha,\beta)$ is sufficiently close to $(1/2,3/2)$, we always have
\bea\label{index equ}
i_{1}(\xi_{\alpha,\beta,\mathfrak{e}})=0,\;\nu_{1}(\xi_{1/2,3/2,\mathfrak{e}})=0;\;i_{-1}(\xi_{\alpha,\beta,\mathfrak{e}})=2,\;\nu_{-1}(\xi_{1/2,3/2,\mathfrak{e}})=0.
\eea
From the Bott iteration formula in Proposition \ref{Bott}, we know
\bea
i_{-1}(\xi_{\alpha,\beta,\mathfrak{e}})=i_{1}(\xi_{\alpha,\beta,\mathfrak{e}})+S_{M}^{+}(1)+
\sum_{0<\theta<\pi}(S_{M}^{+}(e^{\sqrt{-1}\theta})-S_{M}^{-}(e^{\sqrt{-1}\theta}))\nonumber
-S_{M}^{-}(-1),
\eea
where $M=\xi_{\alpha,\beta,\mathfrak{e}}(2\pi)$. Then combine $(\ref{index equ})$ and the Proposition \ref{split norm}, we have
$$
2=\sum_{0<\theta<\pi}(S_{M}^{+}(e^{\sqrt{-1}\theta})-S_{M}^{-}(e^{\sqrt{-1}\theta})).
$$
Also, Proposition \ref{split norm} tells us that the only possible case is $M\approx R(2\pi-\vartheta_{1})\diamond R(2\pi-\vartheta_{2})$ with some
$\vartheta_{1}, \vartheta_{2}\in(0,\pi)$. This completes the proof.
\end{proof}
\section{The linearized system of the planar restricted $(n+1)$-body problem}\label{Se2}
%Considering $n$ particles with masses $m_1,...,m_n$, let $Q=(q_1,...,q_n)\in (\R^2)^n$ be
%the position vector, and $P=(p_1,...,p_n)\in (\R^2)^n$ be the momentum vector.
Let $Q, P\in \mathbb{R}^2$,
from equation (\ref{R3B-E}),
the Hamiltonian function of the massless body has the form
\bea
H(P,Q,t)=\frac{||P||^2}{2}-\mathcal{V}(Q,t),\ \ \mathcal{V}(Q,t)=\sum_{i=1}^{n}\frac{m_{i}}{||r(t)R(\theta(t))a_{i}-Q||}.\nonumber
\eea
%We denote by $\mathbb{J}_n=diag(J_2,...,J_2)_{2n\times2n} $ and  $M=diag (m_1,m_1,m_2,m_2,...,m_n, m_n)_{2n\times2n}$.
Then the corresponding fundamental solution $\ga$ of periodic solution $(P(t), Q(t))$ is given by
\bea \dot{\ga}(t)=J_{4}D^2H(P(t),Q(t))\ga(t),\,\ \ga(0)=I_{4}. \label{ga}  \eea
where
\bea
D^2H(P(t),Q(t))=\left(\begin{array}{cccc} I_{2} & O_{2} \\
O_{2} & -D^2\mathcal{V}(Q(t),t)
\end{array}\right),\nonumber
\eea
$O_{2}$ is the zero matrix.

In the following, we analysis the linear Hamiltonian sytem (\ref{ga}) directly, by using
the linear symplectic transformations, we re-derive (\ref{EREH}) of the linear Hamiltonian system of ERE in the planar restricted $(n+1)$-body problem, see Theorem \ref{red-ERE}.

For the elliptic relative equilibrium (\ref{qN})
$$
q_{N}(t)=r(t)R(\theta(t))a_{N},
$$
which generated by $a_{N}$, where
$$ r(t)=\frac{\Omega^2/\mu}{1+\mathfrak{e}\cos\theta(t)}, \ \ \ \
-\mu a_{n+1}=\sum_{i=1}^{n}m_{i}\frac{a_{i}-a_{n+1}}{||a_{n+1}-a_{i}||^3}.$$
Consider the linear Hamiltonian system at $(p_{n+1}(t), q_{n+1}(t))$ with $P(t)=\dot{Q}(t), Q(t)=q_{n+1}(t)$.
Based on the relation $r^{2}\dot{\theta}(t)=\Omega\neq0$, we first change variable $t$ to the true anomaly $\theta$.
Without loss of generality, we assume $\Omega>0$ and the initial condition $\theta(0)=0$, then $\theta(T)=2\pi$, where $T$ is the minimal period of $x(t)$.

Now define $\tilde{\gamma}(\theta)=\gamma(t(\theta))$, $\tilde{r}(\theta)=r(t(\theta))=\frac{\Omega^2/\mu}{1+\mathfrak{e}\cos\theta}$,
the derivative corresponding to $\theta$ is denoted by $'$, then direct computation shows that
\bea
\tilde{\gamma}'(\theta)=J_{4}\tilde{B}(\theta)\tilde{\gamma}(\theta),\ \ \tilde{\gamma}(0)=I_{4},\label{timechange}
\eea
where
$$
\tilde{B}(\theta)=\left(\begin{array}{cccc} \frac{\tilde{r}^2}{\Omega}I_{2} & O_{2} \\
O_{2} & -\frac{\tilde{r}^2}{\Omega}D^2\mathcal{V}(Q(t(\theta)),t(\theta))
\end{array}\right).
$$
We can further simplify $\tilde{B}(\theta)$, since for any $z\in\mathbb{R}^{2}$, $\mathcal{V}(\tilde{r}R(\theta)z, t(\theta))=\mathcal{U}(z)/\tilde{r}$, where
$$
\mathcal{U}(z)=\sum_{i=1}^{n}\frac{m_{i}}{||a_{i}-z||},
$$
we have
\bea
\frac{\partial^{2}\mathcal{V}(\tilde{r}R(\theta)z, t(\theta))}{\partial z^2}\big{|}_{z=a_{n+1}}=\frac{D^2\mathcal{U}(a_{n+1})}{\tilde{r}},\label{drU1}
\eea
on the other hand, direct computation shows that
\bea
\frac{\partial^{2}\mathcal{V}(\tilde{r}R(\theta)z, t(\theta))}{\partial z^2}\big{|}_{z=a_{n+1}}=\tilde{r}^2R^{\mathcal{T}}(\theta)D^{2}\mathcal{V}(\tilde{r}R(\theta)a_{n+1},t(\theta))R(\theta),\label{drU2}
\eea
compare equalities (\ref{drU1}) and (\ref{drU2}), for $Q(t(\theta))=\tilde{r}R(\theta)a_{n+1}$, we have
$$
D^{2}\mathcal{V}(Q(t(\theta)),t(\theta))=\frac{1}{\tilde{r}^3}R^{-\mathcal{T}}(\theta)D^2\mathcal{U}(a_{n+1})R^{-1}(\theta),
$$
hence
\bea \label{B1}
\tilde{B}(\theta)=\left(\begin{array}{cccc} \frac{\tilde{r}^2}{\Omega}I_{2} & O_{2} \\
O_{2} & -\frac{\mathfrak{R}^{-T}(\theta)D^2\mathcal{U}(a_{n+1})\mathfrak{R}^{-1}(\theta)}{\tilde{r}\Omega}\nonumber
\end{array}\right).
\eea
In order to get Theorem \ref{red-ERE}, we introduce the following lemma to make further transformation,
\begin{lem}\label{symptran}
Let $S(\theta)$ be a path of symplectic matrices depending on $\theta$, if $\hat{\gamma}(\theta)=S^{-1}(\theta)\tilde{\gamma}(\theta)S(0)$, where
$\tilde{\gamma}(\theta)$ is given by (\ref{timechange}) then
$$
\hat{\gamma}'(\theta)=J\Psi_{S}(\tilde{B})\hat{\gamma}(\theta),\;\hat{\gamma}(0)=I,
$$
where $\Psi_{S}(\tilde{B})=S^{T}(\theta)\tilde{B}(\theta)S(\theta)+JS^{-1}(\theta)S'(\theta)$.
\end{lem}
\begin{proof}
The proof is obtained by direct calculation, and we leave it to the readers.
\end{proof}
Now we take
$$
S(\theta)=\left(\begin{array}{cccc} \frac{\sqrt{\Omega}}{\tilde{r}}A^{-T} & O \\
O & \frac{\tilde{r}}{\sqrt{\Omega}}A
\end{array}\right)\left(\begin{array}{cccc} I & \frac{\tilde{r}'}{\tilde{r}}I \\
O & I
\end{array}\right),
$$
where $A\in GL(\R^{2})$ satisfies
$$
J_{2}A = AJ_{2},\ \ A^TA = I_{2}.
$$
Let $\hat{\gamma}(\theta)=S^{-1}(\theta)\tilde{\gamma}(\theta)S(0)$, combine with Lemma \ref{symptran} and the following relations
$$
\tilde{r}(\theta)=\frac{\Omega^2/\mu}{1+\mathfrak{e}\cos\theta},\ \ \tilde{r}'(\theta)=\frac{\mathfrak{e}\sin\theta}{1+\mathfrak{e}\cos\theta}\tilde{r}(\theta),
$$
we have
\bea
\hat{\gamma}'(\theta)=J\Psi_{S}(\tilde{B})\hat{\gamma}(\theta),\ \ \hat{\gamma}(0)=I_{4}.\nonumber
\eea
where
\bea
\Psi_{S}(\tilde{B})=\left( \begin{array}{cccc} I_{2} & O \\
O & I_{2}-\frac{I_{2}+\frac{1}{\mu}R^{-T}(\theta)A^{-1}D^2\mathcal{U}(a_{n+1})AR^{-1}(\theta)}{1+\mathfrak{e}\cos\theta}\nonumber
\end{array}\right).
\eea
Define $\zeta_{\mathfrak{e}}(\theta)=\left(
                            \begin{array}{cc}
                              R^{-1}(\theta) & O \\
                              O & R^{-1}(\theta)\\
                            \end{array}
                          \right)\hat{\gamma}(\theta)
$ and combine with $A^TA = I_{2}$, then we have
\bea
\zeta'_{\mathfrak{e}}(\theta)=J_{4}\left( \begin{array}{cccc} I_{2} & -J_{2} \\
J_{2} & I_{2}-\frac{I_{2}+\frac{1}{\mu}A^{T}D^2\mathcal{U}(a_{n+1})A}{1+\mathfrak{e}\cos\theta}
\end{array}\right)\zeta_{\mathfrak{e}}(\theta),\ \ \zeta_{\mathfrak{e}}(0)=I_{4}.\nonumber
\eea
Denote $\lambda_{\min}, \lambda_{\max}$ be the minimal and maximal eigenvalues of $\frac{1}{\mu}D^2\mathcal{U}(a_{n+1})$ and $x_{1}, x_{2}$ be the unit orthogonal eigenvectors
corresponding $\lambda_{\max}, \lambda_{\min}$.
Since $J_{2}x_{1}$ is orthogonal to $x_{1}$, we obtain $J_{2}x_{1}=\pm x_{2}$. We take $A=(x_{1},J_{2}x_{1})$, then
$$
\frac{1}{\mu}A^{T}D^2\mathcal{U}(a_{n+1})A=diag(\lambda_{\max},\lambda_{\min})=\hat{\alpha}I_{2}+\hat{\beta}\mathcal{N},
$$
with $\hat{\alpha}=\frac{\lambda_{\min}+\lambda_{\max}}{2}$ and $\hat{\beta}=\frac{\lambda_{\max}-\lambda_{\min}}{2}\geq0$. Moreover, direct computation shows that
$$
\mu^{-1}D^2\mathcal{U}(a_{n+1})=-\frac{1}{\mu}\sum_{i=1}^{n}\frac{m_{i}}{||a_{i}-a_{n+1}||^3}I_{2}+\frac{3}{\mu}
\sum_{i=1}^{n}\frac{m_{i}(a_{i}-a_{n+1})(a_{i}-a_{n+1})^{\mathcal{T}}}{||a_{i}-a_{n+1}||^5}.
$$
Therefore,
$$
2\hat{\alpha}=\lambda_{\min}+\lambda_{\max}=Tr(\mu^{-1}D^2\mathcal{U}(a_{n+1}))=\frac{1}{\mu}\sum_{i=1}^{n}\frac{m_{i}}{||a_{i}-a_{n+1}||^3}>0.
$$
This completes the proof of Theorem \ref{red-ERE}.
\begin{rem} From $\hat{\gamma}(2\pi)=\zeta_{\mathfrak{e}}(2\pi)$ and the periodicity of $S(\theta)$, that is $S(2\pi)=S(0)$, we have $\tilde{\gamma}(2\pi)=S(0)\hat{\gamma}(2\pi)S^{-1}(0)=S(0)\zeta_{\mathfrak{e}}(2\pi)S^{-1}(0)$, hence $\tilde{\gamma}(2\pi)$
 and $\zeta_{\mathfrak{e}}(2\pi)$ have the same spectrum, the stability of restricted $(n+1)$ body problem is also determined by $\zeta_{\mathfrak{e}}(2\pi)$.
\end{rem}
%\textbf{Proof of Theorem \ref{unstable-ERE}}: Let $\lambda_{\min}$ be the smallest eigenvalue of $\mu^{-1}D^2\mathcal{U}(a_{N})$, then the elliptic relative equilibria $q_{N}(t)=r(t)R(\theta(t))a_{N}$ of the planar restricted $(n+1)$-body problem is spectral instability for $\lambda_{\min}<-1$, and hyperbolic for $\lambda_{\min}>0$.
\begin{proof}[\textbf{Proof of Theorem \ref{unstable-ERE} and Theorem \ref{stable-ERE}}]
From Theorem \ref{red-ERE}, we know that the fundamental solution of the ERE in the planar restricted $(n + 1)$-body
problem satisfies,
$$
\zeta'_{\mathfrak{e}}(\theta)=J_{4}\mathfrak{B}_{\hat{\alpha},\hat{\beta},\mathfrak{e}}(\theta)\zeta_{\mathfrak{e}}(\theta),\; \zeta_{\mathfrak{e}}(0)=I_{4},
$$
%with
%\bea B(\theta)=\left( \begin{array}{cccc} I_{2} & -J_{2} \\
%J_{2} & I_{2}-\frac{(1+\alpha) I_{2}+\beta\mathcal{N}}{1+\mathfrak{e}\cos\theta}
%\end{array}\right),\quad \theta\in[0,2\pi],\ \ \mathcal{N}=\left(
%                                                             \begin{array}{cc}
%                                                               1 & 0 \\
%                                                               0 & -1 \\
%                                                             \end{array}
%                                                           \right),
% \label{msf} \eea
with $\hat{\alpha}=\frac{\lambda_{\min}+\lambda_{\max}}{2}>0$, $\hat{\beta}=\frac{\lambda_{\max}-\lambda_{\min}}{2}\geq0$, $\lambda_{\min}$ and $\lambda_{\max}$ are the minimal and maximal eigenvalues of
$
\mu^{-1}D^2\mathcal{U}(a_{n+1}).
$

If $\lambda_{\min}<-1$, then $\hat{\beta}>\hat{\alpha}+1$, hence Theorem \ref{stable-thm}-(I) implies that $\zeta_{\mathfrak{e}}(2\pi)\approx D(\lambda)\diamond M$ with $M\in\Sp(2),1\neq\lambda>0$. If $\lambda_{\min}>\max\{-\frac{\lambda_{\max}}{3+3\mathfrak{e}},-1\}$, then $0\leq\hat{\beta}<\min\{\hat{\alpha}+1,k(\mathfrak{e})\hat{\alpha}\}$, hence Theorem \ref{stable-thm}-(II) implies that $\zeta_{\mathfrak{e}}(2\pi)$ is hyperbolic. This complete the proof of
Theorem \ref{unstable-ERE}.

If $\lambda_{\min}>\max\{-\frac{\lambda_{\max}}{2},-1\}$ with $(\lambda_{\max},\lambda_{\min})$ sufficiently close to $(2,-1)$, direct computation shows that it's just
equivalent to $0\leq\hat{\beta}<\min\{3\hat{\alpha},\hat{\alpha}+1\}$ with $(\hat{\alpha},\hat{\beta})$ is sufficiently close to $(1/2,3/2)$, hence Theorem \ref{stable-thm2} implies that $\zeta_{\mathfrak{e}}(2\pi)$ is linear stable and and
$\zeta_{\mathfrak{e}}(2\pi)\approx R(2\pi-\vartheta_{1})\diamond R(2\pi-\vartheta_{2})$ with some $\vartheta_{1},\vartheta_{2}\in(0,\pi)$. This completes the proof of Theorem \ref{stable-ERE}.

\end{proof}

\section{Applications}
For the elliptic relative equilibria  $q_{n+1}(t)=r(t)R(\theta(t))a_{n+1}$ of the planar restricted $(n+1)$-body problem which produced by the planar regular $(1+\mathfrak{n})$-gon central configuration, where $\mathfrak{n}$ particles, each with mass $m$, positioned at the vertices of a regular polygon, rotating around $m_{0}$ at its center. Here note that
$n=\mathfrak{n}+1, M=\mathfrak{n}m$. Then, from \cite{XLZ}, we have following lemma
\begin{lem}(See Proposition 4.1 \cite{XLZ})\label{lem-gon-1}
Suppose that $a=(a_{1},...,a_{n})$ forms a regular $(1+\mathfrak{n})$-gon configuration, then the minimal eigenvalue $\lambda_{\min}$ and the maximal eigenvalue $\lambda_{\max}$ of $\mu^{-1}D^2\mathcal{U}(a_{n+1})$ are given by
\bea\label{rela-eig}
\lambda_{\max}=\frac{A}{\varpi^2}+\frac{|B|}{\varpi^2},\;\lambda_{\min}=\frac{A}{\varpi^2}-\frac{|B|}{\varpi^2}=\frac{l_{3}(x,\vartheta)}{\varpi^2}-1,
\eea
where $\varpi, A, B, l_{3}(x,\vartheta)$ are given by p.305 and p.311 of \cite{BE}. In particular, $\vartheta=\frac{\pi}{\mathfrak{n}}$, when the massless body lies at $S_3$.
$\vartheta=0$, when the massless body lies at $S_1$ and $S_{2}$.
\end{lem}
\begin{lem}(See Proposition 4.3 in \cite{XLZ})\label{lem-gon-2}
If the massless body lies at $S_1$ or $S_{2}$, we have $l_{3}(x,0)<0$.
\end{lem}
\begin{lem}(See Proposition 4.4 in \cite{XLZ})\label{lem-gon-3}
If the massless body lies at $S_3$, we have $l_{3}(x,\frac{\pi}{\mathfrak{n}})>0$. Moreover, when $\frac{m_{0}}{M}\rightarrow+\infty$, we have
$l_{3}(x,\frac{\pi}{\mathfrak{n}})\rightarrow 0^{+}$ and $\frac{A}{\varpi^{2}}\rightarrow (\frac{1}{2})^{+}$.
\end{lem}
%\textbf{Proof of Theorem \ref{1+n-gon}}: If the massless body lies at $S_1$ or $S_{2}$, the elliptic relative equilibrium of the massless body is spectral unstable for any $\mathfrak{e}\in[0,1)$.
In the following, we can see that Theorem \ref{XLZ-TH2} is a special case of our Theorem \ref{stable-ERE}, even though Theorem \ref{XLZ-TH2} has already been proved by Xie, Liu and Zhou in \cite{XLZ}.
\begin{proof}[\textbf{Proof of Theorem \ref{XLZ-TH2}}]
If the massless body lies at $S_3$, Lemma \ref{lem-gon-3} tells us $l_{3}(x,\frac{\pi}{\mathfrak{n}})>0$, from the equalities (\ref{rela-eig}) in Lemma \ref{lem-gon-1}, one can see
$\lambda_{\min}>-1$, $\frac{A}{\varpi^2}+1>\frac{|B|}{\varpi^2}$.
When $\frac{m_{0}}{M}$ is sufficiently large, Lemma \ref{lem-gon-3} also tells us that
$l_{3}(x,\frac{\pi}{\mathfrak{n}})\rightarrow 0^{+}$ and $\frac{A}{\varpi^{2}}\rightarrow (\frac{1}{2})^{+}$, which implies $(\lambda_{\max},\lambda_{\min})\rightarrow(2,-1)$ and
$\frac{3A}{\varpi^{2}}>\frac{A}{\varpi^{2}}+1>\frac{|B|}{\varpi^2}$ , this inequality is equivalent to
$\lambda_{\min}>-\lambda_{\max}/2$. Overall, we have $\lambda_{\min}>\max\{-\lambda_{\max}/2,-1\}$ and $(\lambda_{\max},\lambda_{\min})\rightarrow(2,-1)$ for $\frac{m_{0}}{M}$ is sufficiently large, it satisfies the condition in Theorem \ref{stable-ERE}, hence for any fixed $\mathfrak{e}\in[0,1)$, the elliptic relative equilibrium of the massless body is linear stable, when $\frac{m_{0}}{M}$ is sufficiently large.
\end{proof}
Also, using our more general Theorem \ref{unstable-ERE}, we can prove Theorem \ref{1+n-gon}, which improved the unstable result of Xie, Liu and Zhou in Theorem \ref{XLZ-TH1}.
\begin{proof}[\textbf{Proof of Theorem \ref{1+n-gon}}]
If the massless body lies at $S_1$ or $S_{2}$, from Lemma \ref{lem-gon-1} and Lemma \ref{lem-gon-2}, we have $\lambda_{\min}<-1$, it satisfies the condition in Theorem \ref{unstable-ERE}, hence the elliptic relative equilibrium of the massless body is spectral unstable for any $\mathfrak{e}\in[0,1)$.
\end{proof}
\section{Appendix}
\subsection{A brief review the Maslov-type index theory}\label{App-index}
In this section, we briefly review the index theory for symplectic path.  The detail  could be  found  in \cite{Lon4}. Let $(\mathbb{R}^{2n}, \Omega)$ be the standard symplectic vector space with coordinates $(x_1,\cdots,x_n, y_1,\cdots,y_n)$, then $\Omega=\sum_{i=1}^n dx_i\wedge dy_i$. Let $J=\left(\begin{array}{cc}0 & -I_{n} \\ I_{n} & 0\end{array}\right)$, where $I_n$ is the identity matrix on $\mathbb{R}^{n}$.
The symplectic group $\Sp(2n)$ is defined by
$$  \Sp(2n)=\{M\in GL(2n, \mathbb{R}) | M^\mathcal{T}JM=J \},   $$
with topology induced from $\mathbb{R}^{4n^2}$.  For $\tau>0$ we are interested in paths in $\Sp(2n)$:
$$
\mathcal{P}_{\tau}(2 n)=\left\{\gamma \in C([0, \tau], \operatorname{Sp}(2 n)) \mid \gamma(0)=I_{2 n}\right\}.
$$
For any $\omega\in\mathbb{U}$(the unite circle in $\mathbb{C}$),
the following real function was introduced in  \cite{Lon4},
$$
D_{\omega}(M)=(-1)^{n-1} \bar{\omega}^{n} \operatorname{det}\left(M-\omega I_{2 n}\right).
$$
Thus for any $\omega\in\mathbb{U}$
 the following  $\omega$-degenerate hypersurface of codimension one in $\Sp(2n)$ is defined  \cite{Lon4}:
$$    \Sp(2n)_\omega^0=\{M\in\Sp(2n) | D_{\omega}(M)=0 \}.   $$
Moreover, the $\omega$-regular set of $\Sp(2n)$ is defined by
$$
\begin{gathered}
\operatorname{Sp}(2 n)_{\omega}^{\pm}=\left\{M \in \operatorname{Sp}(2 n) \mid \pm D_{\omega}(M)<0\right\}, \\
\operatorname{Sp}(2 n)_{\omega}^{*}=\operatorname{Sp}(2 n)_{\omega}^{+} \cup \operatorname{Sp}(2 n)_{\omega}^{-}
\end{gathered}.
$$
For any two continuous path $\xi$ and $\eta:[0, \tau] \rightarrow \operatorname{Sp}(2 n)$ with $\xi(\tau)=\eta(0)$, we define their concatenation as:
$$
\eta * \xi(t)= \begin{cases}\xi(2 t), & \text { if } 0 \leqslant t \leqslant \tau / 2 \\ \eta(2 t-\tau), & \text { if } \tau / 2 \leqslant t \leqslant \tau.\end{cases}
$$
%Given any  two $2m_k\times2m_k$ matrices of square block form
%$M_{k}=\left(\begin{array}{ll}A_{k} & B_{k} \\ C_{k} & D_{k}\end{array}\right)$ with $k=1, 2$, the $\diamond$-product of $M_1$ and $M_2$
% is defined by the following $2(m_1+m_2)\times 2(m_1+m_2)$ matrix  $M_{1} \diamond M_{2}$:
%$$
%M_{1} \diamond M_{2}=\left(\begin{array}{cccc}
%A_{1} & 0 & B_{1} & 0 \\
%0 & A_{2} & 0 & B_{2} \\
%C_{1} & 0 & D_{1} & 0 \\
%0 & C_{2} & 0 & D_{2}
%\end{array}\right).
%$$
We define a special continuous symplectic path $\xi_n\subset \Sp(2n)$ by
$$
\xi_{n}(t)=\left(\begin{array}{cc}
2-\frac{t}{\tau} & 0 \\
0 & \left(2-\frac{t}{\tau}\right)^{-1}
\end{array}\right)^{\diamond n} \quad \text { for } 0 \leqslant t \leqslant \tau.
$$

For $M\in \Sp(2n)_\omega^0$, we define a co-orientation of  $\Sp(2n)_\omega^0$
at $M$ by the positive direction $\frac{d}{dt}Me^{tJ}|_{t=0}$ of the path $Me^{tJ}$ with $|t|$ sufficiently small.
Now will give the definition of $\omega$-index.

\begin{defi}\label{def:Maslov-type index} For $\omega\in \mathbb{U}$, $\gamma\in \mathcal{P}_{\tau}$, \textbf{the $\omega$-index} of $\gamma$  can be defined as
$$i_\omega(\gamma)=[e^{-\delta J}(\gamma*\xi_n): \Sp(2n)_\omega^0].$$  For $\delta>0$ small enough, where $[\cdot : \cdot]$ is the intersection number. We also denote $$\nu_\omega(\gamma)=\dim_{\mathbb{C}}\ker_{\mathbb{C}}(\gamma(\tau)-\omega I).$$
and if $\gamma(\tau)\in\Sp(2n)^*_\omega$, we can take $\delta=0$.
\end{defi}
\begin{rem}\label{homtopy}
According to the homotopy invariance of the Maslov index, one have $i_\omega(\gamma)=i_\omega(\gamma^{\epsilon})$, for sufficiently small $\epsilon>0$, where
$\gamma^{\epsilon}(t)=e^{-\epsilon t J}\gamma(t)$.
\end{rem}
%Obviously, If $\gamma(\tau)\in\Sp(2n)^*_\omega$, we can take $\delta=0$. %Since $\xi_n(0)\in\Sp(2n)^+_\omega$, from \cite{Lon4}, we have
%\begin{lem} \label{lem: index even}
%$i_\omega(\gamma)$ is even if and only if  $e^{\vep J}\gamma(\tau)\in\Sp(2n)^+_\omega$ for $\vep>0$ small enough.
%\end{lem}
\begin{prop}\cite[Theorem 8,Page142]{Lon4}\label{per-index}
For any $\tilde{\gamma}$ which are sufficiently $C^{0}$-closed to $\gamma$, we have
$$
|i_\omega(\tilde{\gamma})-i_\omega(\gamma)|\leq \nu_{\omega}(\gamma).
$$
\end{prop}

As in \cite{Lon4}, let $\Omega^{0}(M)$ be the path-connected component containing $M=\gamma(\tau)$ of the set
\bea
\Omega(M)=\{N\in\text{Sp}(2n)\;|\;\sigma(N)\cap\mathbb{U}=\sigma(M)\cap\mathbb{U}\;\text{and}\;\nu_{\omega}(M)=\nu_{\omega}(N),\forall\omega\in\sigma(M)\cap\mathbb{U}\}\nonumber
\eea
where $\sigma(\cdot)$ denotes the spectrum of a matrix, that is the set of its total eigenvalues. Here $\Omega^{0}(M)$ is called the homotopy component of $M$ in $\Sp(2n)$. For a continuous family of paths $\gamma_{s}(t)$ with
$(s,t)\in[0,1]\times[0,\tau]$, $\gamma_{s}(\tau)\in \Omega^{0}(\gamma_{0}(\tau))$, then $i_{\omega}(\gamma_{s})$ is independent of $s$.

In \cite{Lon4}, the following symplectic matrices were introduced as basic normal forms:
	\begin{equation*}
		\begin{aligned}
	&D(\lambda) = \begin{pmatrix} \lambda & 0 \\ 0 & \lambda^{-1} \end{pmatrix}, \quad \lambda = \pm 2;\ \
	N_1(\lambda, a) = \begin{pmatrix} \lambda & a \\ 0 & \lambda \end{pmatrix}, \quad \lambda = \pm 1, a = \pm 1, 0,\\
	&R(\theta) = \begin{pmatrix} \cos\theta & -\sin\theta \\ \sin\theta & \cos\theta \end{pmatrix},\ \
	N_2(\omega, b) = \begin{pmatrix} R(\theta) & b \\ 0 & R(\theta) \end{pmatrix}, \quad \theta \in (0, \pi) \cup (\pi, 2\pi),\\
	    \end{aligned}
	\end{equation*}
	where $ b = \begin{pmatrix} b_1 & b_2 \\ b_3 & b_4 \end{pmatrix} $, $ b_i \in \mathbb{R} $ and $ b_2 \neq b_3 $, $\omega = e^{\sqrt{-1}\theta} $.
\begin{prop}[\cite{Lon4}]\label{homo}
For any $M\in\Sp(2n)$, it can be connected to $N$ in $\Omega^{0}(M)$, where
\bea
N=M_{1}\diamond\cdots\diamond M_{j},\nonumber
\eea
with $M_{i}, i =1,...,j$ in basic normal form, for simplicity, we still denote it as $M\approx N$.
\end{prop}
\begin{prop}[\cite{Lon4}]\label{decom}
For any two paths $\gamma_{j} \in \mathcal{P}_{\tau}(2 n_{j})$ with $j=0, 1$, let $\gamma_{0} \diamond \gamma_{1}(t)=\gamma_{0}(t) \diamond \gamma_{1}(t)$ for all $t\in[0,\tau]$. It is obvious that
\begin{equation*}
\label{eq:IndSum} i_{\omega}\left(\gamma_{1} \diamond \gamma_{2}\right)=i_{\omega}\left(\gamma_{1}\right)+i_{\omega}\left(\gamma_{2}\right), \quad \forall \omega \in \mathbb{U}.
\end{equation*}
\end{prop}
For any $M\in\Sp(2n)$ and $\omega\in\mathbb{U}$, the splitting
number $S^\pm_M(\omega)$ of $M$ at $\omega$ are defined by
\begin{equation}
\label{eq: SplitNumber} S_{M}^{\pm}(\omega)=\lim _{\vep \rightarrow 0^{+}} i_{\omega e^{\pm \sqrt{-1} \vep}}(\gamma)-i_{\omega}(\gamma),
\end{equation}
for any path  $\gamma \in \mathcal{P}_{\tau}(2 n)$ satisfying $\gamma(\tau)=M$.  The splitting number is well defined, i.e. they are independent of the choice of path $\gamma$.  From the property of Maslov-type index, if $M_0, M_1\in\Sp(2n)$ which are symplectic similar, i.e. $M_0\approx M_1$, then
$$S_{M_{0}}^\pm(\omega)= S_{M_{1}}^\pm(\omega). $$
For any $M_i\in\Sp(2n_i)$ with $i=0$ and $1$, there holds
$$
S_{M_{0} \diamond M_{1}}^{\pm}(\omega)=S_{M_{0}}^{\pm}(\omega)+S_{M_{1}}^{\pm}(\omega), \quad \forall \omega \in \mathbb{U}.
$$

We only list the splitting number of the normal form of $\Sp(2)$, which is enough for the current paper.
For $M\in\Sp(2)$, it is symplectic similar to one of the following  normal forms:
\begin{equation}\label{definition of R(theta)}
\begin{array}{c}
D(\lambda)=\left(\begin{array}{cc}
\lambda & 0 \\
0 & \lambda^{-1}
\end{array}\right), \quad \lambda\in \mathbb{R}\setminus\{0,\pm1\} \\
N_{1}(\lambda, a)=\left(\begin{array}{cc}
\lambda & a \\
0 & \lambda
\end{array}\right), \quad \lambda=\pm 1, a=\pm 1,0 \\
R(\theta)=\left(\begin{array}{cc}
\cos \theta & -\sin \theta \\
\sin \theta & \cos \theta
\end{array}\right), \quad \theta \in(0, \pi) \cup(\pi, 2 \pi)
\end{array}
\end{equation}

\begin{prop}\label{split norm} (\cite{Lon4}) For $M\in\Sp(2n)$ and $\omega\in\mathbb{U}$, $0<\theta<\pi$, there hold
\begin{equation}
\begin{aligned}
S_{M}^{\pm}(\omega) &=0, \quad \text { if } \omega \notin \sigma(M), \\
S_{M}^{\pm}(\omega) &=S_{M}^{\mp}(\bar{\omega}), \\
0 \leqslant S_{M}^{\pm}(\omega) & \leqslant \nu_\omega(M), \\
S_{M}^{+}(\omega)+S_{M}^{-}(\omega) & \leqslant \operatorname{dim} \operatorname{ker}(M-\omega I)^{2 n}, \quad \text { if } \omega \in \sigma(M), \\
\left(S_{N_{1}(1, a)}^{+}(1), S_{N_{1}(1, a)}^{-}(1)\right) &=\left\{\begin{array}{ll}
(1,1), & \text { if } a=0,1, \\
(0,0), & \text { if } a=-1,
\end{array}\right.\\
\left(S_{N_{1}(-1, a)}^{+}(-1), S_{N_{1}(-1, a)}^{-}(-1)\right) &=\left\{\begin{array}{ll}
(1,1), & \text { if } a=0,-1, \\
(0,0), & \text { if } a=1,
\end{array}\right.\\
\left(S_{R(\theta)}^{+}\left(e^{\sqrt{-1} \theta}\right), S_{R(\theta)}^{-}\left(e^{\sqrt{-1} \theta}\right)\right) &=(0,1) \\
\left(S_{R(2 \pi-\theta)}^{+}\left(e^{\sqrt{-1} \theta}\right), S_{R(2 \pi-\theta)}^{-}\left(e^{\sqrt{-1} \theta}\right)\right) &=(1,0),
\end{aligned}
\end{equation}
\end{prop}
We also have following useful Bott iteration formula,
\begin{prop}(\cite[$P_{204}$]{Lon4})\label{Bott}
For any symplectic path $\gamma\in\mathcal{P}_{\tau}(2n)$, $M=\gamma(\tau)$, and $\omega=e^{i\vartheta}$ with $\vartheta\in[0,2\pi)$, we have
\bea
i_{\omega}(\ga)=i_{1}(\ga)+S_{M}^{+}(1)+
\sum_{0<\theta<\vartheta}(S_{M}^{+}(e^{\sqrt{-1}\theta})-S_{M}^{-}(e^{\sqrt{-1}\theta}))\nonumber
-S_{M}^{-}(\omega),
\eea
\end{prop}
%For any symplectic path $\gamma\in\mathcal{P}_{\tau}(2n)$ and $m\in\mathbb{Z}^+$, we define the $m$-th iteration $\gamma^m: [0,m\tau]\to\Sp(2n)$ by
%$$
%\gamma^{m}(t)=\gamma(t-j \tau) \gamma(\tau)^{j}, \quad \text { for } j \tau \leqslant t \leqslant(j+1) \tau, j=0,1, \ldots, m-1.
%$$
%
%\begin{defi} \label{def:mean index} The \textbf{mean index} of $\gamma$ is defined by
%\bea \hat{i}(\gamma)=\lim_{m\to\infty}\frac{i_1(\gamma^m)}{m}. \eea
%
%\end{defi}
%
%Meanwhile the Bott-type iteration formula shows
%$$
%i_{z}\left(\gamma^{m}\right)=\sum_{\omega^{m}=z} i_{\omega}(\gamma), \;\; \forall z\in\mathbb{U},
%$$
%
%
%
%which implies
%$$
%\hat{i}(\gamma) =\frac{1}{2 \pi} \int_{0}^{2 \pi} i_{e^{\sqrt{-1} \theta}}(\gamma) d \theta.
%$$
%
%Y. Long has showed a more explicit form of the mean index by the Bott-type formula and the basic norm form of $\Sp(2n)$. We give this form as follows.
%\begin{lem}\label{lem: mean index form}(\cite{Lon4}) For any symplectic path $\gamma\in\mathcal{P}_{\tau}(2n)$ and $M=\gamma(\tau)$, we have
%\begin{equation}
%\hat{i}(\gamma)=i(\gamma)+S_{M}^{+}(1)-C(M)+\sum_{\theta \in(0,2 \pi)} \frac{\theta}{\pi} S_{M}^{-}\left(e^{\sqrt{-1} \theta}\right),
%\end{equation}
%where  $C(M)=\sum_{0<\theta<2 \pi} S_{M}^{-}\left(e^{\sqrt{-1} \theta}\right)$.

%\end{lem}
\subsection{The relationship with the relative Morse index}\label{App-index-relation}
The relative Morse index could be defined by spectral flow \cite{HS09}. The spectral flow was introduced by Atiyah, Patodi and
Singer \cite{APS} in their study of index theory on manifolds with
boundary.
Let $\{A(\theta),\theta\in[0,1]\}$ be a continuous path of
self-adjoint Fredholm operators on a Hilbert space $\mathcal{H}$. Roughly
speaking, the spectral flow of path $\{A(\theta),\theta\in[0,1]\}$
counts the net change in the number of negative eigenvalues of
$A(\theta)$ as $\theta$ goes from $0$ to $1$, where the enumeration
follows from the rule that each negative eigenvalue crossing to the
positive axis contributes $+1$ and each positive eigenvalue crossing
to the negative axis contributes $-1$, and for each crossing, the
multiplicity of eigenvalue is counted. More precisely, as shown in
\cite{APS},  let
$$\wp=\bigcup_{\theta\in[0,1]} \sigma (A({\theta})),
$$ where $\sigma(A(\theta))$ is the spectrum for $A(\theta)$, then $\wp$ is a closed subset of the $(\theta,\lambda)$-plane. The spectral flow $Sf(\{A({\theta})\})$
 is defined to be the intersection number of $\wp$ with the line $\lambda=-\epsilon$ with respect
 to the usual orientation for some small positive $\epsilon$.

Now, consider the following Hamiltonian system
\begin{eqnarray} \dot{z}(t)= JB(t)z(t), \quad z(0)= z(T).  \label{linea system}
\end{eqnarray} Denote by $A=-J\frac{d}{dt},$ which is densely
defined in the Hilbert space $E=L^2([0,T];\mathbb C^{2n})$ with the
domain
\begin{eqnarray*}
D=\left\{\left.z(t)\in W^{1,2}([0,T];\mathbb C^{2n})\,\right|\,z(0)=z(T)\right\}.
\end{eqnarray*}
 $B$ is a bounded linear operator defined by $(Bz)(t)=B(t)z(t)$ on $E$. Then $A$ is a self-adjoint Fredholm operator with compact resolvent.
Then the relative Morse index $I(A,A-B)$ is define by
\begin{eqnarray*}
I(A,A-B)=-Sf(\{A-sB, s\in[0,1]\}).
\end{eqnarray*}
%We come back to the Hamiltonian systems, suppose $B(s,t)\in C([0,1]\times[0,T], S(2n))$. For $s\in[0,1]$,  let $B_s\in\mathcal{B}(2n)$. For such two operators $A-B_0$ and $A-B_1$, we can define the relative Morse index via spectral flow.
%In fact, by \cite{HS}, we have,
%\begin{eqnarray*}
%I(A-B_0,A-B_1)=-Sf(\{A-B(s), s\in[0,1]\}). \label{5a.8}
%\end{eqnarray*}
\begin{prop}\cite[Theorem 2.5 and Lemma 4.6]{HS09}\label{index relation}
For  imaginary number $\nu$  such that   $\omega=e^{\nu T}\in \mathbb U\setminus \{1\}$, we have
\begin{eqnarray*} I(A, A-B)=i_1(\gamma)+n.  \label{5a.18}  \end{eqnarray*}
and
\begin{eqnarray*} I(A-\nu J, A-\nu J-B)=i_\omega(\gamma),  \label{5a.20}  \end{eqnarray*}
where $\gamma$ is the fundamental solution of system (\ref{linea system}) with $\gamma(0)=I$.
\end{prop}

\medskip

%\noindent {\bf Acknowledgements.} The author thanks Professor Penghui Wang for allowing her to freely choose her research topic. The author also thanks Professor Xijun Hu's team for the seminar they offered, in which she learned about index theory and the $N$-body problem.

\vspace{2mm}
\noindent {\textbf{Declaration:}}
About a year ago, the author noticed that the study of the linear stability of the elliptic relative equilibria for massless particle under the regular $(1+\mathfrak{n})$-gon
central configuration in paper \cite{XLZ} (First version: https://arxiv.org/abs/2310.00286), specifically Theorem \ref{XLZ-TH1}, had not been completely resolved. This sparked the author's interest. Through parameter analysis, the author found that this system can be embedded into a more general $(\alpha,\beta)$-system, which was first introduced in paper \cite{HLO}. Within this system, the problem becomes much clearer. In particular, the author realized that when fixing $\alpha$ or $\beta$, or taking $\beta=k\alpha$, the index of the system possesses monotonicity, which is a highly favorable property (Proposition \ref{monot} and \ref{monot2}). Most importantly, the author found that when $(\alpha,\beta)$ lies in certain specific regions, the stability analysis of this system shares similarities with the analyses of the Lagrangian solutions and Euler solutions in papers \cite{HLS} and \cite{ZL}, respectively. Also, in the critical case $\beta=\alpha+1$, the monodormy matrix of the system relate to a two by two system which is more clearer (Lemma \ref{cri-eig}). Therefore, through further analysis of the existing results, we obtain the stability results for the $(\alpha,\beta)$-system and apply them to the restricted $(n+1)$-body problem. In particular, the final part completely resolves the instability problem in Theorem \ref{XLZ-TH1}, accomplishing the initial motivation of our work.

In view of the recent impact of AI on mathematics, the author must declare that all the research began a year ago and
all the theoretical research and results in this article were obtained by the author through analysis based on existing literature, without relying on any AI.
The figure was drawn using MATLAB. The author takes responsibility for this declaration.

%no AI was used in the theoretical research and results of this problem; they were obtained by the author through analysis based on existing literature.
%All results in this article were obtained by the author after a study of index theory and the results of Paper \cite{HW} and \cite{XLZ}, without relying on any AI. The author takes responsibility for this declaration.

\end{document}